\documentclass[12pt,twoside,leqno]{article}
\usepackage{amsmath}
\usepackage{amssymb}
\usepackage{amsxtra}
\usepackage{amscd}
\usepackage{amsthm}
\usepackage[mathscr]{eucal}

\theoremstyle{plain}

\newtheorem{sbthm}[subsubsection]{Theorem}
\newtheorem{sbprop}[subsubsection]{Proposition}
\newtheorem{sbcor}[subsubsection]{Corollary}

\theoremstyle{definition}

\newtheorem{sbpara}[subsubsection]{}

\newenvironment{pf}{\proof[\proofname]}{\endproof}

\begin{document}

\title{Extended period domains, algebraic groups, and higher Albanese manifolds}

\author
{Kazuya Kato
\footnote{
Department of mathematics, University of Chicago, 
Chicago, Illinois, 60637, USA.},
Chikara Nakayama
\footnote{
Department of Economics, Hitotsubashi University, 
2-1 Naka, Kunitachi, Tokyo 186-8601, Japan.},
Sampei Usui
\footnote{
Graduate School of Science, Osaka University,
Toyonaka, Osaka, 560-0043, Japan.}}

\maketitle
\renewcommand{\mathbb}{\bold}

\newcommand\Cal{\mathcal}
\newcommand\define{\newcommand}
\define\gp{\mathrm{gp}}%
\define\fs{\mathrm{fs}}%
\define\an{\mathrm{an}}%
\define\mult{\mathrm{mult}}%
\define\add{\mathrm{add}}%
\define\Ker{\mathrm{Ker}\,}%
\define\Coker{\mathrm{Coker}\,}%
\define\Hom{\mathrm{Hom}\,}%
\define\Ext{\mathrm{Ext}\,}%
\define\rank{\mathrm{rank}\,}%
\define\gr{\mathrm{gr}}%
\define\cHom{\Cal{Hom}}
\define\cExt{\Cal Ext\,}%

\define\cC{\Cal C}
\define\cD{\Cal D}
\define\cO{\Cal O}
\define\cS{\Cal S}
\define\cM{\Cal M}
\define\cG{\Cal G}
\define\cH{\Cal H}
\define\cE{\Cal E}
\define\cF{\Cal F}
\define\cN{\Cal N}
\define\fF{\frak F}
\define\fg{\frak g}
\define\fh{\frak h}
\define\Dc{\check{D}}
\define\Ec{\check{E}}

\newcommand{\N}{{\mathbb{N}}}
\newcommand{\Q}{{\mathbb{Q}}}
\newcommand{\Z}{{\mathbb{Z}}}
\newcommand{\R}{{\mathbb{R}}}
\newcommand{\C}{{\mathbb{C}}}
\newcommand{\bN}{{\mathbb{N}}}
\newcommand{\bQ}{{\mathbb{Q}}}
\newcommand{\bF}{{\mathbb{F}}}
\newcommand{\bZ}{{\mathbb{Z}}}
\newcommand{\bP}{{\mathbb{P}}}
\newcommand{\bR}{{\mathbb{R}}}
\newcommand{\bC}{{\mathbb{C}}}
\newcommand{\bS}{{\bold{S}}}
\newcommand{\bbQ}{{\bar \mathbb{Q}}}
\newcommand{\ol}[1]{\overline{#1}}
\newcommand{\too}{\longrightarrow}
\newcommand{\respect}{\rightsquigarrow}
\newcommand{\compatible}{\leftrightsquigarrow}
\newcommand{\upc}[1]{\overset {\lower 0.3ex \hbox{${\;}_{\circ}$}}{#1}}
\newcommand{\Gmlog}{\bG_{m, \log}}
\newcommand{\Gm}{\bG_m}
\newcommand{\ep}{\varepsilon}
\newcommand{\Spec}{\operatorname{Spec}}
\newcommand{\val}{{\mathrm{val}}} 
\newcommand{\n}{\operatorname{naive}}
\newcommand{\bs}{\operatorname{\backslash}}
\newcommand{\Gal}{\operatorname{{Gal}}}
\newcommand{\gal}{{\rm {Gal}}({\bar \Q}/{\Q})}
\newcommand{\galp}{{\rm {Gal}}({\bar \Q}_p/{\Q}_p)}
\newcommand{\gall}{{\rm{Gal}}({\bar \Q}_\ell/\Q_\ell)}
\newcommand{\wep}{W({\bar \Q}_p/\Q_p)}
\newcommand{\wel}{W({\bar \Q}_\ell/\Q_\ell)}
\newcommand{\Ad}{{\rm{Ad}}}
\newcommand{\BS}{{\rm {BS}}}
\newcommand{\even}{\operatorname{even}}
\newcommand{\End}{{\rm {End}}}
\newcommand{\odd}{\operatorname{odd}}
\newcommand{\GL}{\operatorname{GL}}

\newcommand{\np}{\text{non-$p$}}
\newcommand{\g}{{\gamma}}
\newcommand{\G}{{\Gamma}}
\newcommand{\Lam}{{\Lambda}}
\newcommand{\La}{{\Lambda}}
\newcommand{\lam}{{\lambda}}
\newcommand{\la}{{\lambda}}
\newcommand{\uL}{{{\hat {L}}^{\rm {ur}}}}
\newcommand{\uQp}{{{\hat \Q}_p}^{\text{ur}}}
\newcommand{\sel}{\operatorname{Sel}}
\newcommand{\dt}{{\rm{Det}}}
\newcommand{\Sig}{\Sigma}
\newcommand{\fil}{{\rm{fil}}}
\newcommand{\SL}{{\rm{SL}}}
\newcommand{\spl}{{\rm{spl}}}
\newcommand{\st}{{\rm{st}}}
\newcommand{\Isom}{{\rm {Isom}}}
\newcommand{\Mor}{{\rm {Mor}}}
\newcommand{\bg}{\bar{g}}
\newcommand{\id}{{\rm {id}}}
\newcommand{\cone}{{\rm {cone}}}
\newcommand{\al}{a}
\newcommand{\ChL}{{\cal{C}}(\La)}
\newcommand{\Image}{{\rm {Image}}}
\newcommand{\toric}{{\operatorname{toric}}}
\newcommand{\torus}{{\operatorname{torus}}}
\newcommand{\Aut}{{\rm {Aut}}}
\newcommand{\Qp}{{\mathbb{Q}}_p}
\newcommand{\barQp}{{\mathbb{Q}}_p}
\newcommand{\Qpur}{{\mathbb{Q}}_p^{\rm {ur}}}
\newcommand{\Zp}{{\mathbb{Z}}_p}
\newcommand{\Zl}{{\mathbb{Z}}_l}
\newcommand{\Ql}{{\mathbb{Q}}_l}
\newcommand{\Qlur}{{\mathbb{Q}}_l^{\rm {ur}}}
\newcommand{\F}{{\mathbb{F}}}
\newcommand{\eps}{{\epsilon}}
\newcommand{\epsLa}{{\epsilon}_{\La}}
\newcommand{\epsLaVxi}{{\epsilon}_{\La}(V, \xi)}
\newcommand{\epsOLaVxi}{{\epsilon}_{0,\La}(V, \xi)}
\newcommand{\Qplin}{{\mathbb{Q}}_p(\mu_{l^{\infty}})}
\newcommand{\otimesQplin}{\otimes_{\Qp}{\mathbb{Q}}_p(\mu_{l^{\infty}})}
\newcommand{\galFl}{{\rm{Gal}}({\bar {\Bbb F}}_\ell/{\Bbb F}_\ell)}
\newcommand{\gallur}{{\rm{Gal}}({\bar \Q}_\ell/\Q_\ell^{\rm {ur}})}
\newcommand{\galFF}{{\rm {Gal}}(F_{\infty}/F)}
\newcommand{\galFv}{{\rm {Gal}}(\bar{F}_v/F_v)}
\newcommand{\galF}{{\rm {Gal}}(\bar{F}/F)}
\newcommand{\epsVxi}{{\epsilon}(V, \xi)}
\newcommand{\epsOVxi}{{\epsilon}_0(V, \xi)}
\newcommand{\plim}{\lim_
{\scriptstyle 
\longleftarrow \atop \scriptstyle n}}
\newcommand{\sig}{{\sigma}}
\newcommand{\ga}{{\gamma}}
\newcommand{\del}{{\delta}}
\newcommand{\Vss}{V^{\rm {ss}}}
\newcommand{\Bst}{B_{\rm {st}}}
\newcommand{\Dpst}{D_{\rm {pst}}}
\newcommand{\Dcrys}{D_{\rm {crys}}}
\newcommand{\DdR}{D_{\rm {dR}}}
\newcommand{\Fin}{F_{\infty}}
\newcommand{\Kla}{K_{\lambda}}
\newcommand{\Ola}{O_{\lambda}}
\newcommand{\Mla}{M_{\lambda}}
\newcommand{\Det}{{\rm{Det}}}
\newcommand{\Sym}{{\rm{Sym}}}
\newcommand{\LaSa}{{\La_{S^*}}}
\newcommand{\cX}{{\cal {X}}}
\newcommand{\MHG}{{\frak {M}}_H(G)}
\newcommand{\tauMla}{\tau(M_{\lambda})}
\newcommand{\Fvur}{{F_v^{\rm {ur}}}}
\newcommand{\Lie}{{\rm {Lie}}}
\newcommand{\LMH}{{\rm {LMH}}}
\newcommand{\cB}{{\cal {B}}}
\newcommand{\cL}{{\cal {L}}}
\newcommand{\cW}{{\cal {W}}}
\newcommand{\fq}{{\frak {q}}}
\newcommand{\cont}{{\rm {cont}}}
\newcommand{\SC}{{SC}}
\newcommand{\Om}{{\Omega}}
\newcommand{\dR}{{\rm {dR}}}
\newcommand{\crys}{{\rm {crys}}}
\newcommand{\hatSig}{{\hat{\Sigma}}}
\newcommand{\rdet}{{{\rm {det}}}}
\newcommand{\ord}{{{\rm {ord}}}}
\newcommand{\Alb}{{{\rm {Alb}}}}
\newcommand{\BdR}{{B_{\rm {dR}}}}
\newcommand{\BdRO}{{B^0_{\rm {dR}}}}
\newcommand{\Bcrys}{{B_{\rm {crys}}}}
\newcommand{\Qw}{{\mathbb{Q}}_w}
\newcommand{\barkappa}{{\bar{\kappa}}}
\newcommand{\cP}{{\Cal {P}}}
\newcommand{\cZ}{{\Cal {Z}}}
\newcommand{\oppLa}{{\Lambda^{\circ}}}
\newcommand{\bG}{{\mathbb{G}}}
\newcommand{\br}{{{\bold r}}}
\newcommand{\triv}{{\rm{triv}}}
\newcommand{\sub}{{\subset}}
\newcommand{\LD}{{D^{\star,\mild}_{\SL(2)}}}
\newcommand{\LbD}{{D^{\star}_{\SL(2)}}}
\newcommand{\dbDv}{{D^{\star}_{\SL(2),\val}}}
\newcommand{\nspl}{{{\rm nspl}}}
\newcommand{\lval}{{[\val]}}
\newcommand{\mild}{{{\rm{mild}}}}
\newcommand{\lan}{\langle}
\newcommand{\ran}{\rangle}
\newcommand{\rar}{{{\rm {rar}}}}
\newcommand{\rac}{{{\rm {rac}}}}
\newcommand{\Rep}{{\operatorname{Rep}}}

\begin{center} Dedicated to Professor Steven Zucker on his 65th birthday. 
\end{center}
\medskip

\begin{abstract} 
 For a linear algebraic group $G$ over $\Q$, we consider the period domains $D$ classifying $G$-mixed Hodge structures, and construct the extended 
period domains 
$D_{\Sig}$. 
We give an interpretation of higher Albanese manifolds by Hain and Zucker by using the above $D$ for some $G$, and extend them via $D_\Sigma$.

\end{abstract}

\renewcommand{\thefootnote}{\fnsymbol{footnote}}
\footnote[0]{Primary 14C30; 
Secondary 14D07, 32G20} 

\setcounter{section}{-1}
\section{Introduction}\label{s:intro}

For a linear algebraic group $G$ over $\Q$, we consider the period domains $D$ for  $G$-mixed Hodge structures. We construct the extended 
period domains 
$D_{\Sig}$, 
the space of nilpotent orbits. 
In this paper, we give an interpretation of higher Albanese manifolds by Hain and Zucker by using the above $D$ for some $G$, and extend them via $D_\Sigma$.

In Section \ref{s:hz}, we review a work on higher Albanese manifolds by Hain and Zucker in \cite{HZ1}.
In Section \ref{s:D}, we define $D$ by modifying the definition of Shimura variety over $\C$ by Deligne (\cite{De}). 
  In Section \ref{s:DSig}, we introduce the extended period domain $D_{\Sig}$ and state the main results \ref{t:property}, \ref{t:main}. 
  In Section \ref{s:ex}, we explain the relation of this $D_{\Sig}$ with the theory for 
the usual period domains (\cite{KNU}), and 
as examples, the Mumford--Tate domains, and 
mixed Shimura varieties.
  In the case where $D$ is the pure Mumford--Tate domain (cf.\ Green--Griffiths--Kerr's book \cite{GGK}), 
$D_{\Sigma}$ essentially coincides with the one by Kerr--Pearlstein (\cite{KP}). 
In Section \ref{s:hab}, we study higher Albanese manifolds by Hain and Zucker by using the present $D$ for some $G$, and extend them via $D_\Sigma$.
  
  We omit the details of constructions and proofs of the general theory in Section \ref{s:DSig} 
  in this paper, which are to be published elsewhere. 
  
 For the $p$-adic variant of this paper, Koshikawa and the first author are preparing \cite{KK}.

\section{Relation of the work \cite{HZ1} and the present article}\label{s:hz}

In this section, we briefly review the work of Hain and Zucker on unipotent variations of mixed Hodge structure \cite{HZ1} with which we compare our present article.

\subsection{The main theorem of \cite{HZ1}}

Let $X$ be a connected smooth algebraic variety over $\C$ and $b\in X$, and let $\Bbb F$ be a real field, the field of rational numbers,  or the ring of integers.
Then, \cite{HZ1} (1.6) asserts that
there is an equivalence of categories:
$$
\begin{pmatrix}
\text{good unipotent variations of mixed}\\
\text{Hodge structure on $X$ with}\\
\text{unipotency $\le r$, defined over $\Bbb F$}
\end{pmatrix}
\overset{\sim}{\to}
\begin{pmatrix}
\text{mixed Hodge theoretic}\\
\text{representations of $\C\pi_1(X,b)/J^{r+1}$,}\\
\text{defined over $\Bbb F$}
\end{pmatrix}
$$
where $J$ is the augmentation ideal, i.e., the kernel of $\varepsilon:\C\pi_1(X,b)\to\C$, $\gamma \mapsto 1$ ($\gamma\in \pi_1(X,b)$).

An outline of the proof is as follows.
\begin{sbpara}\label{lr} 
The functor from the left-hand-side to the right-hand-side is given by taking the monodromy representation on the fiber over the base point $b$.
\end{sbpara}

\begin{sbpara}\label{rl} 
The correspondence from the right-hand-side to the left-hand-side is given by using higher Albanese manifold of $X$.
\end{sbpara}

\begin{sbpara}\label{rigid} 
The rigidity of variations of mixed Hodge structure is shown under \lq\lq good" condition (\cite{HZ1} (1.5)) at the boundary of $X$.
This rigidity ensures that \ref{rl} yields the inverse functor of \ref{lr}.
\end{sbpara}


\subsection{Iterated integration theory of Chen \cite{Ch}}\label{s:chen}

We review the result of Chen in \cite{Ch}. 

\begin{sbpara}\label{dfls} 
Let $I$ be the interval $[0,1]$. 
A loop on $X$ with base point $b$ is a $C^{\infty}$ map $\gamma:I\to X$ with $\gamma(0)=\gamma(1)=b$.
Let $PX= P_bX$ be the loop space on $(X,b)$ which is the topological space consisting of all loops on $X$ with base point $b$ endowed with compact-open topology.

A local parameter system of $PX$ is a pair $(U,\phi)$ of an open set $U$ of $\bR^n$ and a map $\phi:U\to PX$ such that $\varphi:I\times U\to X$ with $\varphi(t,u):=\phi(u)(t)$ is a $C^{\infty}$ map.
Here $n$ is a non-negative integer.
For an open set $V$ of $\R^m$ and a $C^{\infty}$ map $f:V\to U$, $(V,\phi\circ f)$ is also a local parameter system.
For a non-negative integer $k$, a $k$-differential form on $ PX$ is a collection $\omega=(\omega_\phi)_{\phi}$ of a $k$-differential form $\omega_\phi$ on $U$ for a local parameter system $(U,\phi)$ such that, for $f:U\to V$ as above, $f^*\omega_\phi=\omega_{\phi\circ f}$. 
\end{sbpara}

\begin{sbpara}\label{itin} 
Let $r$ be a positive integer.
Let $\Delta_r=\{(t_1,\dots,t_r)\in\R^r\,|\, 0\le t_1\le\dots\le t_r\le1\}$ be the $r$-simplex.
Let $\pi_j:X^r\to X$ be the $j$-th projection, $1\le j\le r$.
Let $\pi:\Delta_r\times PX\to PX$ be the projection and let $\varphi:\Delta_r\times PX\to X^r$ be the map defined by $\varphi(t_1,\dots,t_r,\gamma):=(\gamma(t_1),\dots,\gamma(t_r))$.
Let $A^k(X)$ and $A^k( PX)$ be $k$-th forms on $X$ and on $ PX$, respectively.
For positive integers $p_{1},\dots,p_{r}$, put $q=\sum_{j=1}^{r}(p_j-1)$.
The iterated integral 
$$
I:A^{p_1}(X)\times\cdots\times A^{p_r}(X)\to A^{q}( PX),\quad
(\omega_1,\dots,\omega_r)\mapsto \int \omega_1\dots\omega_r,
$$
along $\CD  PX @<\pi<<\Delta_r\times PX@>\varphi>>X^r\endCD$, is defined as follows.

It is enough to define it on each local parameter system $\phi:U\to  PX$ with $U\subset\R^n$.
Let $\iota_j$ be the contraction of a differential form with $\frac{\partial}{\partial t_j}$ and set 
$\alpha_j:=\iota_j\varphi^*\pi_j^*\omega_j$, $1\le j\le r$.
Then $\alpha_j$ is a $(p_j-1)$-th form on $\Delta_r\times U$.
Write
$$
\alpha_1\wedge\dots\wedge\alpha_r
=:\sum_{1\le i_1\le\dots\le i_q\le n}\alpha_{i_1\dots i_q}(t_1,\dots,t_r,u_1,\dots,u_n)du_{i_1}\wedge \dots\wedge du_{i_q}.
$$
Define
$$
\int \omega_1\dots\omega_r
:=\sum_{1\le i_1\le\dots\le i_q\le n}(\int_{\Delta_r}\alpha_{i_1\dots i_q}dt_1\dots dt_r)
du_{i_1}\wedge \dots\wedge du_{i_q}.
$$
\end{sbpara}

\begin{sbpara}\label{bcom}

The bar complex $B^\bullet(X)$ is the subcomplex of the de Rham complex $A^\bullet(PX)$ on $PX$, which is determined by the de Rham complex $A^\bullet(X)$ of $X$ as follows.

$B^{q}(X)$ is the subspace of $A^{q}(PX)$ generated by the images of 
$$
I:A^{p_1}(X)\times\cdots\times A^{p_r}(X)\to A^{q}( PX),\quad
(\omega_1,\dots,\omega_r)\mapsto \int \omega_1\dots\omega_r,
$$
for all positive integers $r$ and $p_{1},\dots,p_{r}$ such that $q=\sum_{j=1}^{r}(p_{j}-1)$.

For $1\le j\le r$, let $\nu_{j}:=\sum_{k=1}^{j}(p_{k}-1)$, and let $\nu_0:=0$.
The exterior differential $d:B^q(X)\to B^{q+1}(X)$ is described as  
$$
d\int \omega_1\dots\omega_r
=\sum_{j=1}^{r}(-1)^{\nu_{j-1}+1}\int\omega_{1}\dots\omega_{j-1}d\omega_{j}\omega_{j+1}\dots\omega_{r}
$$
$$
+\sum_{j=1}^{r-1}(-1)^{\nu_{j}+1}\int\omega_{1}\dots\omega_{j-1}(\omega_{j}\wedge\omega_{j+1})\omega_{j+2}\dots\omega_{r}.
$$

Length filtration $(L_{r}B^{\bullet}(X))_{r}$ is an increasing filtration by subcomplexes of $B^{\bullet}(X)$, where $L_{r}B^{\bullet}(X)$ is generated by 
$\int\omega_{1}\dots\omega_{j}$ for differential forms $\omega_{1},\dots,\omega_{j}$ on $X$ of positive degree over $j\le r$.

Denote by $B^{\bullet}(X)_\C$ etc.\ the $\C$-valued iterated integrals etc.
Hodge filtration $(F^{p}B^{\bullet}(X)_\C)_{p}$ is a decreasing filtration by subcomplexes of $B^{\bullet}(X)_\C$, where $F^{p}B^{\bullet}(X)_\C$ is generated by 
$\int\omega_{1}\dots\omega_{r}$ for differential forms $\omega_{1}\in F^{p_{1}}A^{\bullet}(X)_\C,\dots,\omega_{r}\in F^{p_{r}}A^{\bullet}(X)_\C$ of positive degree such that $\sum_{j}p_{j}\ge p$.
\end{sbpara}

\begin{sbpara}\label{bcoh}
An element of $B^{0}(X)$ is a function on the loop space $PX$:
$$
B^0(X)\times PX\to \R, \quad
(\int \omega_1\dots\omega_r, \gamma)\mapsto \int_\gamma \omega_1\dots\omega_r.
$$
In this case, every $\omega_{j}$ is a 1-form on $X$ and, by writing $\gamma^{*}\omega_{j}=:f_{j}(t)dt$, it is described as
$$
\int_\gamma \omega_1\dots\omega_r
=\int_{0\le t_{1}\le\dots\le t_{r}\le1}f_{1}(t_{1})\dots f_{r}(t_{r})
dt_{1}\dots dt_{r}.
$$

This induces
$$
H^0(B^\bullet(X))\times \pi_1(X,b)\to\R
$$
and
$$
H^0(B^\bullet(X))\to\Hom(\Z\pi_1(X,b),\R)
$$
and also
$$
L_rH^0(B^\bullet(X))\to\Hom(\Z\pi_1(X,b)/J^{r+1},\R).
$$

The theorem of Chen in \cite{Ch} asserts that the last homomorphism is an isomorphism.

The filtrations on $\C \pi_1(X,b)/J^{r+1}$, which are induced from the length filtration and the Hodge filtration on the iterated integrals $B^\bullet(X)$, form a mixed Hodge structure called the $r$-th {\it canonical variation of mixed Hodge structure} as $b$ varies over $X$ (\cite{HZ1} (4.21)).

\end{sbpara}


\subsection{Higher Albanese manifolds in \cite{HZ1}}\label{hamHZ}

Put $G:=\pi_{1}(X,b)$.
Let $R$ be a commutative ring with unity.
Let $\varepsilon:RG\to R$ be the augmentation map.
Let $\Delta:RG\to RG\otimes RG$ be the coproduct $\Delta(g):=g\otimes g$.
Let $RG^{\wedge}:=\underset{\leftarrow\;r}{\lim}RG/J^{r+1}$ be the $J$-adic completion and $\hat J$ the closed ideal of $RG^{\wedge}$ generated by $J$.
Define
$$
\hat G_{R}:=\{h\in RG^{\wedge}\,|\,
\varepsilon(h)=1,\, \Delta(h)=h\hat\otimes h\}
\subset 1+\hat J,
$$
$$
\frak g_{R}:=\{h\in\hat J\,|\,\Delta(h)=1\hat\otimes h+h\hat\otimes 1\}.
$$
Let $\hat G_{r,R}:=\hat G_R/(\hat G_R\cap(1+\hat J^{r+1}))$ and
$\frak g_{r,R}:=\frak g_R/(\frak g_R\cap\hat J^{r+1})$ (\cite{HZ1} (2.13)).
Let $F$ be the Hodge filtration on $\frak g_{r,\C}$ induced by the one on $\C G/J^{r+1}$.
Let $F^{0}\hat G_{r,\C}$ be the corresponding subgroup of $\hat G_{r,\C}$.
The higher Albanese manifold in \cite{HZ1} (5.15) is defined by
$$
\Alb^{r}(X):=\hat G_{r,\Z}\bs\hat G_{r,\C}/F^{0}\hat G_{r,\C}.
$$


\subsection{Inverse correspondence}\label{inv}

We explain \ref{rl}.
An inverse correspondence is constructed as follows (\cite{HZ1} (5.21)).
Given a mixed Hodge theoretic representation $V$, i.e., a ring homomorphism $\C \pi_1(X,b)/J^{r+1}\to\End(V)$ of mixed Hodge structures, we have a map $\Alb^r(X)\to D(V)$, where $D(V)$ is the classifying space of Hodge filtrations on $V$.
Composing the higher Albanese map $X\to\Alb^r(X)$ with the above map and pulling back the universal variation of mixed Hodge structure on the classifying space, we get a variation of mixed Hodge structure on $X$.

\subsection{The aim of the present paper}\label{aim}
In this paper, we have the following two contributions \ref{result1} 
and \ref{result2} to the work of Hain and Zucker.

\begin{sbpara}\label{result1} 
We give a description of the functor 
represented by the higher Albanese manifold in terms of tensor 
functors (see Section \ref{s:hab1}). 
Here we give a rough sketch of it (the precise 
and more general statement is given in Theorem \ref{thm2}). 
 
Let $\Gamma_r$ be the image of $\pi_1(X,b)\to \C\pi_1(X, b)/J^{r+1}$. 
By \cite{P} p.85, p.474, cf.\  also \cite{HZ1} (2.17) iii), if we define the subgroups $\Gamma^r$ of $\pi_1(X, b)$ by 
$\Gamma^0:= \pi_1(X, b)$ and $\Gamma^{i+1}:= [\pi_1(X,b), \Gamma^i]$ 
for $i\geq 0$, then $\Gamma_r$ is the quotient group of 
$\pi_1(X,b)/\Gamma^r$ such that the kernel of $\pi_1(X, b)/\Gamma^r 
\to \Gamma_r$ consists of all elements of $\pi_1(X, b)/\Gamma^r$ of 
finite orders. 
Let $\cC_{X,\Gamma_{r}}$ be the category of variations of $\Q$-MHS $\cH$ on $X$ 
satisfying the following conditions. 
 
(i) For any $w\in \Z$, 
$\gr^W_w\cH$ is a constant  polarizable Hodge structure. 
 
(ii) $\cH$ is good at infinity in the sense of \cite{HZ1} (1.5). 
 
(iii) The monodromy action of $\pi_1(X, b)$ on $\cH_{\Q}(b)$ (which is 
unipotent under (i)) factors through $\Gamma_r$. 
 
Then our result is roughly that for a complex analytic space $S$, 
there is a functorial bijection between the set $\text{Mor}(S, 
\Alb^r(X))$ of morphisms and the set of isomorphism classes of exact 
tensor functors $\cC_{X,\Gamma_{r}}\to \text{MHS}(S)$, where $\text{MHS}(S)$ denotes the category of analytic families of $\Q$-MHS 
parametrized by $S$, which sends $h_X$ to $h_S$ for any $\Q$-MHS $h$ (more precisely, see \ref{moduli1} (i)).
 Here objects 
of $\text{MHS}(S)$ need not satisfy Griffiths transversality, though objects of 
$\cC_{X,\Gamma_r}$ should satisfy it. 
$h_X$ (resp.\ $h_S$) denotes the 
constant variation (resp.\ family) of $\Q$-MHS on $X$ (resp.\ $S$) associated to 
$h$.

\end{sbpara} 
 
\begin{sbpara}\label{result2} 
We construct toroidal partial 
compactifications of $\Alb^r(X)$, and describe the functors 
represented by them generalizing \ref{result1} to its log version. 
See Section \ref{s:tpc} for details. 
\end{sbpara} 
 
\begin{sbpara} 
We will deduce these results \ref{result1} and \ref{result2} from the 
work of Hain and Zucker and from our general theory of (extended) period 
domains for mixed Hodge structures associated to algebraic groups, 
which we develop in Section \ref{s:D} and Section \ref{s:DSig}. 
 
\end{sbpara}

\section{The period domain $D$}\label{s:D}
  Let $G$ be a linear algebraic group over $\Q$.
  Let $G_u$ be the unipotent radical of $G$. 
  Let $\mathrm{Rep}(G)$ be the category of 
finite-dimensional linear representations of $G$ over $\Q$. 

\subsection{$G$-mixed Hodge structures and period domain $D$}
  We define and consider the notion of {\it $G$-mixed Hodge structure} ({\it $G$-MHS} for short).

\begin{sbpara}
  As in \cite{De}, let $S_{\C/\R}$ be the Weil restriction of the multiplicative group ${\mathbb G}_m$ from $\C$ to $\R$. It represents the functor $R\mapsto (\C\otimes_\R R)^\times$ for commutative rings $R$ over $\R$.
  We have  $S_{\C/\R}(\R)=\C^\times$, and $S_{\C/\R}$ is understood as $\C^\times$ regarded as an algebraic group over $\R$.

 Let $w: {\mathbb G}_{m,\R}\to S_{\C/\R}$ be the homomorphism induced from the natural maps $R^\times\to (\C\otimes_\R R)^\times$ 
 for commutative rings $R$ over $\R$.
\end{sbpara}

\begin{sbpara}\label{SCR} 
A linear representation of $S_{\C/\R}$ over $\R$ is equivalent to a finite-dimensional $\R$-vector space $V$ endowed with a decomposition
$$V_\C:=\C\otimes_\R V=\bigoplus_{p,q\in \Z} \; V_\C^{p,q}$$
such that for any $p,q$, $V_\C^{q,p}$ coincides with the complex conjugate of $V_{\C}^{p,q}$ (that is, the image of $V_\C^{p,q}$ under $\C\otimes_\R V\to \C\otimes_\R V\;;\;z\otimes v\mapsto {\bar z}\otimes v$). For a linear representation $V$ of $S_{\C/\R}$, the corresponding decomposition is defined by 
$$V_\C^{p,q}=\{v\in V_\C\;|\; [z]v= z^p{\bar z}^qv\;\text{for}\;z\in \C^\times\}.$$
Here $[z]$ denotes $z$ regarded as an element of $S_{\C/\R}(\R)$. 
\end{sbpara}

\begin{sbpara}
\label{D} 
  Let $h_0: S_{\C/\R}\to (G/G_u)_\R$ be a homomorphism.
  Assume that the composite ${\mathbb G}_{m,\R}\overset{w}\to S_{\C/\R}\to (G/G_u)_\R$ is $\Q$-rational and central. 
  Assume also that for one (and hence any) lifting ${\mathbb G}_{m,\R} \to G_{\R}$ of this composite, 
the adjoint action of ${\mathbb G}_{m,\R}$ on $\Lie(G_u)_\R=\bR \otimes_{\bQ}\Lie(G_u)$ is of weight $\leq -1$. 

  Then, for any $V \in \Rep(G)$, the action of ${\mathbb G}_m$  on $V$ via a 
lifting ${\mathbb G}_m \to G$ of the above ${\mathbb G}_m\to G/G_u$  defines 
a rational increasing filtration $W$ on $V$ called the {\it weight filtration}, which is independent of the lifting. 

In the above situation, a {\it $G$-mixed Hodge structure} ({\it $G$-MHS}, for short) is defined as an exact $\otimes$-functor from $\mathrm{Rep}(G)$ 
to the category of $\Q$-MHS keeping the underlying vector spaces with the weight filtrations. 
\end{sbpara}

\begin{sbpara}\label{hpure}
Let $H$ be a $G$-MHS. 
  By \ref{SCR} and Tannaka duality (\cite{DM}), the Hodge decompositions of $\gr^W$ of $H(V)$ for $V\in \Rep(G)$ 
give a homomorphism $S_{\C/\R}\to (G/G_u)_\R$ such that the composite ${\mathbb G}_{m,\R}\overset{w}\to S_{\C/\R}\to (G/G_u)_\R$ is $\Q$-rational and central. 
  We call this homomorphism the {\it homomorphism associated with $H$}. 
\end{sbpara}

\begin{sbpara}\label{pd}
  We define the {\it period domain $D$ associated to $G$ and $h_0$} as the set of all isomorphism classes of $G$-MHS whose associated homomorphism $S_{\C/\R}\to (G/G_u)_\R$ is $(G/G_u)(\R)$-conjugate to $h_0$. 
  In this case, $D$ is also called the {\it period domain of type $(G,h_0)$.}
\end{sbpara}

\subsection{Complex analytic structure of $D$}

We first give a real analytic understanding of $D$, and then consider the complex analytic structure of it.

Fix a homomorphism $h_0: S_{\C/\R}\to (G/G_u)_\R$ as in \ref{D}.

\begin{sbpara} Let $h: S_{\C/\R}\to G_\R$ be a homomorphism such that the composite 
${\mathbb G}_{m,\R}\overset{w}\to S_{\C/\R}\overset h \to G_\R \to (G/G_u)_\R$ coincides with $h_0 \circ w$. 
We define an $\R$-subspace $L(h)$ of $\Lie(G)_\R=\R\otimes_\Q \Lie(G)$ as the set of all $\delta\in \Lie(G)_\R$ such that 
the $(p,q)$-Hodge component of $\delta$ (\ref{SCR}) with respect to the adjoint action of $S_{\C/\R}$ via $h$ is $0$ unless $p<0$ and $q<0$. 
\end{sbpara}

\begin{sbpara}\label{Rprop}  
%
%
%
%
For $\delta\in L(h)$, we obtain a $G$-MHS $H(h,\delta)$ as follows. For a linear representation $V$ of $G$ over $\R$, $H(h,\delta)(V)$ is $(V, W, F)$, where $W$ is 
the weight filtration on $V$ (\ref{D}) and $F$ is the Hodge filtration on $V_\C=\C\otimes_\R V$ defined in the following way. 
Let 
$V_\C=\bigoplus_{p,q} V_\C^{p,q}$ be the Hodge decomposition  (\ref{SCR}) defined by the action of $S_{\C/\R}$ via $h$.  
Let $$F^p :=\exp(i\delta)\bigoplus_{p'\geq p,\,q\in \Z}V_\C^{p',q}.$$
\end{sbpara}

\begin{sbprop}\label{hdelta}  
  The above construction $(h, \delta)\mapsto H(h, \delta)$ gives a bijection from the set of all $(h, \delta)$ as above onto the set of all isomorphism classes of $G$-MHS. 
\end{sbprop}

\begin{pf} By \cite{CKS} and by Tannaka duality (cf.\ \cite{DM}). 
\end{pf}

\begin{sbpara}
Consider the action of the subgroup $G(\R)G_u(\C)$ of $G(\C)$ on $D$ defined by changing Hodge filtrations.
\end{sbpara}
\begin{sbprop}\label{tst} The action of 
$G(\R)G_u(\C)$ on $D$ is transitive. 
\end{sbprop}

\begin{pf} This follows from the definition of $D$ in \ref{pd} and 
Proposition \ref{hdelta}.
\end{pf}

\begin{sbpara} \label{dcheck}
Let $\cC$ be the category of triples $(V, W, F)$, where $V$ is a finite-dimensional $\Q$-vector space, $W$ is an increasing filtration on $V$ 
(called the weight filtration), and $F$ is a decreasing filtration on $V_\C$ (called Hodge filtration).  

Let $Y$ be  the set of all isomorphism classes of exact $\otimes$-functors from $\Rep(G)$ to the category $\cC$ preserving the underlying vector spaces and the weight 
filtrations. 

Then $G(\C)$ acts on $Y$ by changing the Hodge filtration. We have $D\subset Y$ and $D$ is stable in $Y$ under the action of $G(\R)G_u(\C)$. 

Let $$\Dc:=G(\C)D \subset Y.$$ Since the action of $G(\C)$ on $\Dc$ is transitive and the isotropy group of each point of $\Dc$ is an algebraic subgroup of $G(\C)$, $\Dc$ has a natural structure of a complex analytic manifold as a $G(\bC)$-homogeneous space.
\end{sbpara}

\begin{sbprop}\label{open}
$D$ is open in $\Dc$. 
\end{sbprop}

\begin{pf} 
  Let $G_r=G/G_u$. 
  By considering the Hodge decomposition of $\Lie(G_r)_\C$, we can see the equality 
$\Lie(G_r)_{\C} = \Lie(G_r)_\R+ F^0\Lie(G_r)_\C$.
Since $\Lie(G_{r})_{\R}\cap F^{0}\Lie(G_{r})_{\C}=0$ in $\Lie(G_{r})_{\C}$, we have $\dim_{\R}\Lie(G_{r})_{\C}/F^{0}\Lie(G_{r})_{\C}=\dim_{\R}\Lie(G_{r})_{\R}$, and hence the proposition 
follows.
\end{pf}

\begin{sbcor}
  $D$ is a complex analytic manifold. 
\end{sbcor}

\subsection{Polarizability} 
  For a linear algebraic group $G$, let $G'$ be the commutator algebraic subgroup.
  
\begin{sbpara}\label{pol} 
  Let $h_0: S_{\C/\R}\to (G/G_u)_\R$ be as in \ref{D}. 
Let $C$ be the image of $i\in \C^\times = S_{\C/\R}(\R)$ by $h_0$ in $(G/G_u)(\R)$. 
  We say that $h_0$ is {\it $\R$-polarizable} if 
$\{a\in (G/G_u)'(\R)\;|\; Ca=aC\}$ is a maximal compact subgroup of $(G/G_u)'(\R)$.

\end{sbpara}

\begin{sbpara}\label{pol2}
  A relationship with the usual $\bR$-polarizability is as follows (\cite{De2} 2.11). 
  Let $h_0$ be as in \ref{D}. 
  Let $H$ be a $G$-MHS such that the associated $S_{\C/\R}\to (G/G_u)_\R$ is $\R$-polarizable. 
  Let $V \in \Rep(G)$. 
  Then for each $w\in \Z$, there is an $\bR$-bilinear form on $\gr^W_w(V)_\R$ which is stable under $(G/G_u)'$ and which polarizes $\gr^W_wH(V)$. 
     
\end{sbpara}

\begin{sbpara}\label{Gamma} 
We will often consider a subgroup $\Gamma$ of $G(\Q)$ satisfying the following condition.

There is a faithful $V\in \Rep(G)$ and a $\Z$-lattice $L$ in $V$ such that $L$ is stable under the action of $\Gamma$.

\end{sbpara}

\begin{sbprop}
  Let $h_0: S_{\C/\R}\to (G/G_u)_\R$ be as in $\ref{D}$.
  Assume that $h_0: S_{\C/\R}\to (G/G_u)_\R$ is $\R$-polarizable ($\ref{pol}$).
  Let $\Gamma$ be a subgroup of $G(\Q)$ satisfying the condition in $\ref{Gamma}$.

Then the following holds. 

$(1)$ The action of $\Gamma$ on $D$ is proper, and the quotient space $\Gamma \bs D$ is Hausdorff. 

$(2)$ If $\Gamma$ is torsion-free and if $\gamma p=p$ with $\gamma\in \Gamma$ and  some $p\in D$, then $\gamma=1$.

$(3)$ If $\Gamma$ is torsion-free, then the projection $D\to \Gamma \bs D$ is a local homeomorphism.
\end{sbprop}

\begin{pf}
(1) 
  By the assumption of $\bold R$-polarizability, 
the action of $\Gamma$ on $D$ is proper.
  By \cite{KNU} Part III 4.2.4.1, the quotient space $\Gamma \bs D$ is Hausdorff. 

(2) By the condition in $\ref{Gamma}$, $\Gamma$ is discrete.
(2) follows then from $\R$-polarizability and torsion-freeness of $\Gamma$.

(3) Since $\Gamma$ is discrete and $D$ is Hausdorff, (2) implies (3) by \cite{KNU} Part III 4.2.4.2.
\end{pf}

\section{Space of nilpotent orbits $D_{\Sig}$}\label{s:DSig}
We define the extended period domain $D_\Sig\supset D$ as the space of nilpotent orbits,  
and state the main results.
  We fix $G$ and $h_0$ as in \ref{D}. 
  Assume that $h_0$ is $\bR$-polarizable.

\subsection{Definition of $D_{\Sig}$}

\begin{sbpara}\label{nilp1}

  A {\it nilpotent cone} is a subset $\sig$ of $\Lie(G)_\R$ 
  satisfying the following (i)--(iii).
  
  (i) $\sig=\R_{\geq 0}N_1+\dots +\R_{\geq 0}N_n$ for some $N_1,\dots,N_n\in \Lie(G)_\R$. 
  
  (ii) For any $V\in \Rep(G)$, the image of $\sig$ under the induced map $\Lie(G)_\R\to \End_\R(V)$ consists of  nilpotent operators. 
  
  (iii) $[N,N']=0$ for any $N, N'\in \sig$. 
  \end{sbpara}
  
  \begin{sbpara}\label{nilp2}
    Let $F\in \Dc$ and let $\sig$ be a nilpotent cone.
     We say that the pair $(\sig, F)$ {\it generates a nilpotent orbit} if the following (i)--(iii) are satisfied.
     
     (i)  There is a faithful $V\in \Rep(G)$ such that the action of $\sig$ on $V_{\bR}$ is admissible with respect to $W$, i.e., 
     there exist a family $(M(\tau,W))_\tau$ of finite increasing filtrations $M(\tau,W)$ on $V$ given for each face $\tau$ of $\sigma$ which satisfy the  compatibility conditions (1)--(4) in \cite{KNU} Part III 1.2.2.

     (ii) $N F^p\subset F^{p-1}$ for any $N\in \sig$ and $p\in \Z$.
     
     (iii) Let $N_1,\dots, N_n$ be as in (i) in \ref{nilp1}. Then $\exp(\sum_{j=1}^n z_jN_j)F\in D$ if $z_j\in \C$ and $\text{Im}(z_j)\gg 0$ ($1\leq j\leq n$). 
     
  A {\it nilpotent orbit} is a pair 
$(\sig, Z)$ of a nilpotent cone $\sigma$ and an $\exp(\sig_\C)$-orbit $Z$ in $\Dc$
satisfying that for any $F\in Z$, $(\sig, F)$ generates a nilpotent orbit. Here $\sig_\C$ denotes the $\C$-linear span of $\sig$ in $\Lie(G)_\C$. 
  
  \end{sbpara}

\begin{sbpara}\label{fan}

  A {\it weak fan $\Sig$ in $\Lie(G)$} is a nonempty set of sharp rational nilpotent cones 
satisfying the conditions that it is closed under taking faces 
and that any $\sig, \sig' \in \Sig$ coincide if they have a common interior 
point and if there is an $F\in \Dc$ such that both $(\sig, F)$ and $(\sig',F)$ generate nilpotent 
orbits. 

For a weak fan $\Sigma$, let $D_{\Sig}$ be the set of all nilpotent orbits $(\sig, Z)$ such that $\sig\in \Sig$.
  Then $D$ is naturally embedded in $D_{\Sig}$ by $F\mapsto(\{0\},F)$.

  Let $\Gamma$ be a subgroup of $G(\Q)$ satisfying the condition in \ref{Gamma}. 
  We say that $\Sigma$ and $\Gamma$ are {\it strongly compatible} if $\Sig$ is stable under 
the adjoint action of $\Gamma$ and if any $\sig \in \Sig$ is generated by elements whose $\exp$ in $G(\R)$ belong to $\Gamma$. 
  If this is the case, $\Gamma$ naturally acts on $D_{\Sig}$. 
\end{sbpara}

\subsection{Log mixed Hodge structures}

\begin{sbpara}\label{glmhl}
We work in the category $\cB(\log)$ of locally ringed spaces over 
$\C$ with fs log structures satisfying a certain condition, which 
contains the category of fs log analytic spaces over $\C$ (\cite{KNU} Part III 1.1). 
For an object $S=(S,\cO_{S},M)$ of $\cB(\log)$, there exists the associated ringed space $S^{\log}=(S^{\log},\cO_{S}^{\log})$ and a proper surjective morphism $S^{\log}\to S$ of ringed spaces (\cite{KNU} Part III).
We denote by $\LMH (S)$ the category of log $\Q$-mixed Hodge structures over $S$ (\cite{KNU} Part III 1.3). 
  
  Let $\Gamma$ be a subgroup of $G(\Q)$ satisfying the condition in \ref{Gamma}. 
A $G$-LMH over $S$ with a {\it $\Gamma$-level structure} is a pair $(H,\mu)$ of an exact $\otimes$-functor $H:\Rep(G)\to\LMH (S)$ and a global section $\mu$ 
of the quotient sheaf $\Gamma\bs {\cal I}$, where $\cal I$ is the following sheaf on $S^{\log}$. For an open set $U$ of $S^{\log}$, ${\cal I}(U)$ is the set of all isomorphisms $H_\Q|_U\overset{\cong}\to \text{id}$ of $\otimes$-functors from $\Rep(G)$  to the category of local systems of $\Q$-modules over $U$ preserving the weight filtrations. 
       
     \end{sbpara}
     \begin{sbpara}\label{type}
     Let $(G, h_0)$ be as in \ref{D}, let $\Gamma$ be a subgroup of $G(\Q)$ satisfying the condition in \ref{Gamma} and let 
  $\Sig$ be a weak fan in $\Lie(G)$ which is strongly compatible with $\Gamma$. 
A $G$-LMH over $S$ with a $\Gamma$-level structure $(H,\mu)$
 is said to be {\it of type $(h_0, \Sig)$} if the following (i) and (ii) are satisfied for any $s\in S$ and any $t\in s^{\log}$.
Take a $\otimes$-isomorphism $\tilde \mu_t: H_{\Q,t}\cong  \;\text{id}$ which belongs to $\mu_t$.

  (i) There is a $\sig\in \Sig$ such that the logarithm of the action of the local monodromy cone 
$\Hom((M_S/\cO^\times_S)_s, \N)\subset \pi_1(s^{\log})$ on $H_{\Q,t}$ is contained, via $\tilde \mu_t$, in $\sigma \subset \Lie(G)_\R$.

  (ii) Let $\sig\in \Sig$ be the smallest cone satisfying (i). Let  $a: \cO_{S,t}^{\log}\to \C$ be a ring homomorphism which induces the evaluation $\cO_{S,s}\to \C$ at $s$ and consider the element $F: V\mapsto {\tilde \mu}_t (a(H(V)))$ of $Y$ (\ref{dcheck}). Then this element belongs to $\Dc$ and $(\sig, F)$ generates a nilpotent orbit (\ref{nilp2}).
  
  If $(H, \mu)$ is of type $(h_0, \Sig)$, we have a map $S \to \Gamma \bs D_{\Sig}$, called the {\it period map} associated to $(H, \mu)$,  which sends $s\in S$ to the class of the nilpotent orbit $(\sig, Z)\in D_{\Sig}$ obtained in the above (ii). 
  
  \end{sbpara}

\begin{sbpara}
Let $(G,h_0, \Gamma, \Sig)$ be as in \ref{type}. We endow $\Gamma \bs D_{\Sig}$ with a topology, a sheaf of rings $\cO$ over $\C$ and a log structure $M$ defined as follows. 
The topology  is the strongest topology for which the period map $S \to \Gamma \bs D_{\Sig}$ is continuous for any $(S, H, \mu)$, where $S$ is an object of $\cB(\log)$, $H$ is a $G$-LMH on $S$, and $\mu$ is a $\Gamma$-level structure which is of  type $(h_0,\Sig)$. For an open set $U$
of $\Gamma \bs D_{\Sig}$, $\cO(U)$ (resp.\ $M(U)$) is the set of all $\C$-valued functions $f$ on $U$ such that for any $(S, H, \mu)$ as above with the period map $\varphi: S \to \Gamma \bs D_{\Sig}$, the pullback of $f$ on $U':=\varphi^{-1}(U)$ belongs to the image of $\cO_{U'}$ (resp.\ $M_{U'}$) in the sheaf of $\C$-valued functions on $U'$.

These structures of $\Gamma \bs D_{\Sig}$ are defined also by defining spaces $E_{\sig}$ ($\sig\in \Sig$) in a similar way as  \cite{KNU} Part III. We get the same structures when we use only $S=E_{\sig}$ for $\sig \in \Sig$ and the universal objects $(H, \mu)$ over $E_{\sig}$ in the above definitions of the structures.
\end{sbpara}

\begin{sbpara}\label{GMHS} Let  $S$ be an object of $\cB(\log)$. Let $S^{\circ}$ be the underlying locally ringed space over $\C$ of $S$ with the trivial log structure. 

By an MHS over $S$, we mean an LMH over $S^{\circ}$.

  Let $(G, h_0)$ be as in \ref{D} and let $\Gamma$ be a subgroup of $G(\Q)$ satisfying the condition in \ref{Gamma}. By a $G$-MHS over $S$ with $\Gamma$-level structure, we mean a 
  $G$-LMH over $S^{\circ}$ with $\Gamma$-level structure. By a {\it $G$-MHS over $S$ with $\Gamma$-level structure of type $h_0$}, we mean a $G$-LMH over $S^{\circ}$ with $\Gamma$-level structure of type $(h_0, \Sig)$ where $\Sig$ is the fan consisting of the one cone $\{0\}$. 

\end{sbpara} 

\subsection{Main results}
We state main results for moduli of $G$-log mixed Hodge structures in our general theory.

\begin{sbthm}
\label{t:property}
Let $(G, h_0, \Gamma, \Sig)$ be as in \ref{type}. Assume that $h_0$ is $\R$-polarizable. 
Then

$(1)$ $\Gamma \bs D_{\Sig}$ is Hausdorff.

$(2)$ When $\Gamma$ is neat, $\Gamma \bs D_{\Sig}$ is a log manifold ({\rm \cite{KNU} Part III 1.1.5}). 
In particular, $\Gamma 
\bs D_{\Sig}$ belongs to $\cB(\log)$.
\end{sbthm}

Here we say that $\Gamma$ is {\it neat} if there is a faithful $V\in \Rep(G)$ such that for any $\gamma\in \Gamma$, the subgroup of $\C^\times$ generated by all eigenvalues of $\gamma: V_{\C}\to V_{\C}$ is torsion-free. 

\begin{sbpara} 
  The outline of the proof is as follows. 
  As in \cite{KNU}, we can define various spaces $D_{\SL(2)}$, $D_{\BS}$, $E_{\sig}$ etc., and have the theory of CKS map. 
  Then, as in \cite{KNU}, 
by using the CKS map, good properties of
$\Gamma \bs D_{\Sig}$ are deduced from those of the space of $\SL(2)$-orbits $D_{\SL(2)}$, which  
reduce to 
the $\bR$-polarizable version of \cite{KNU}. 
  We remark that what were shown in \cite{KNU} by using $\bQ$-polarizations still hold 
under $\bR$-polarizations (\ref{pol2}) .

\end{sbpara}

\begin{sbthm}\label{t:main}
Let $(G, h_0, \Gamma, \Sig)$ be as in Theorem \ref{t:property}.
When $\Gamma$ is neat, 
$\Gamma \bs D_{\Sig}$ represents the contravariant functor from $\cB(\log)$ to {\rm (Set):}

$S\mapsto \{$isomorphism class of 
$G$-LMH over $S$ with a $\Gamma$-level structure of type $(h_0,\Sigma)$ $\}$.
\end{sbthm} 

The proof of \ref{t:main} is similar to the proof of \cite{KNU} Part III, 2.6.6.

\medskip

Concerning extensions of period maps to the boundary, we have: 

\begin{sbthm}\label{t:extend}
Let $(G,h_{0})$ be as in \ref{D}.
Assume that $h_0$ is $\R$-polarizable. 
Let $S$ be a connected, log smooth, fs log analytic space, and let $U$ be the open subspace of $S$ consisting of all points of $S$ at which the log structure of $S$ is trivial.
Let $\G$ be a subgroup of $G(\Q)$ as in \ref{Gamma}.
Assume that $\Gamma$ is neat.

Let $(H,\mu)$ be a $G$-MHS over $U$ with a $\Gamma$-level structure 
of  type $h_0$ (\ref{GMHS}).
Let $\varphi: U\to \G\bs D$ be the associated period map. 
Assume that $(H,\mu)$ extends to a $G$-LMH over $S$ with a $\Gamma$-level structure (\ref{glmhl}). 
Then{\rm:}

\medskip

{\rm(1)} For any point $s\in S$, there exist an open neighborhood
$V$ of $s$, a log modification $V'$ of $V$ $(\text{\cite{KU}}\ 3.6.12)$, a subgroup $\G'$ of $\G$, and a fan $\Sig$ in $\Lie(G)$
which is strongly compatible with $\G'$ such that the period map
$\varphi|_{U\cap V}$ lifts to a morphism $U\cap V \to \G'\bs D$ which
extends uniquely to a morphism $V'\to\G'\bs D_\Sig$ of log manifolds.
Furthermore, we can take a 
commutative group $\Gamma'$.
$$
\CD
U&\supset&U\cap V&\quad\subset\quad& V'\\
@V{\varphi}VV@VVV@VVV\\
\G \bs D@<<< \G' \bs D&\subset& \G' \bs D_\Sig.
\endCD
$$

{\rm(2)} Assume $S\smallsetminus U$ is a smooth divisor. 
Then we can take $V=V'=S$ and $\G'=\G$. 
That is, we have a commutative diagram 
$$
\CD
U&\subset&S\\
@V{\varphi}VV@VVV\\
\G \bs D&\quad\subset\quad& \G \bs D_\Sig.
\endCD
$$

{\rm(3)} Assume that $\G$ is commutative.
Then we can take $\G'=\G$.

\medskip

{\rm (4)} 
Assume that $\G$ is commutative and that the following condition {\rm(i)} is satisfied.
\smallskip

{\rm(i)} There is a finite family $(S_j)_{1\leq j\leq n}$ of connected
locally closed analytic subspaces of $S$ such that $S=\bigcup_{j=1}^n S_j$
as a set and such that, for each $j$, the inverse image of the sheaf
$M_S/\cO^\times_S$ on $S_j$ is locally constant.
\smallskip

Then we can take $\G'=\G$ and $V=S$. 
\end{sbthm} 

Note that in (1), we can take a fan $\Sig$ (we do not need a weak fan). 
\medskip

This  is the $G$-mixed Hodge version of \cite{KNU} Part III, Theorem 7.5.1 in mixed case and of \cite{KU} Theorem 4.3.1 in pure case. 
The proof goes exactly in the same way as in the pure case treated in \cite{KU}. 

\medskip

\section{Basic examples}
\label{s:ex}
We discuss basic examples of $D$ to which our theory can be applied so that we can give $D_{\Sig}$ for these $D$. 

\subsection{Usual period domains}
  We explain that the classical Griffiths domains \cite{Gr} and their mixed Hodge generalization in \cite{U} are 
essentially regarded as 
special cases of the period domains of this paper. 
  In this case, our partial compactifications essentially coincide with those in \cite{KNU} Part III.

  Let $\Lambda=(H_0, W, (\langle \;,\;\rangle_w)_w, (h^{p,q})_{p,q})$ be as usual as in \cite{KNU} Part III. 
  Let $G$ be the subgroup of $\Aut(H_{0,\Q}, W)$ consisting of elements which induce {\it similitudes for $\langle\;,\;\rangle_w$} for each $w$. That is, 
$
G:= \{g\in \Aut(H_{0,\Q}, W)\;|\;$ for any $w$, there is a $t_w\in {\mathbb G}_m$ such that 
$\langle gx,gy\rangle_w = t_w\langle x, y \rangle_w$ for any $x,y \in \gr^W_w\}.
$
  Let $G_1:=\Aut(H_{0,\Q}, W, (\langle \;,\;\rangle_w)_w) \subset G$.

  Let $D(\Lambda)$ be the period domain of \cite{U}. 
Then $D(\Lambda)$ is identified with an open and closed part 
 of  $D$ 
in this paper as follows. 

Assume that $D(\Lam)$ is not empty and fix an $\br\in D(\Lambda)$. Then the Hodge decomposition of  $\gr^W\br$  induces 
$h_0: S_{\C/\R} \to (G/G_u)_\R$. (We have
$\langle [z]x, [z]y\rangle_w= |z|^{2w}\langle x, y\rangle_w$ for $z\in \C^\times$ (see \ref{SCR} for $[z]$).) 
  Consider the associated period domain $D$ (\ref{D}). Then $D$ is a finite disjoint union of $G_1(\R)G_u(\C)$-orbits which are open and closed in $D$. 
  Let $\cD$ be the $G_1(\R)G_u(\C)$-orbit in $D$ consisting of points whose associated homomorphisms $S_{\C/\R}\to (G/G_u)_\R$  are $(G_1/G_u)(\R)$-conjugate to $h_0$. 
  Then the map $H\mapsto H(H_{0,\Q})$ gives a $G_1(\R)G_u(\C)$-equivariant  isomorphism $\cD\overset{\cong}\to D(\Lambda)$. 
  
\subsection{Mixed Mumford--Tate domains}
\label{ss:MT}

\begin{sbpara}\label{MTgp}

  Let $H$ be a $\Q$-MHS whose $\gr^W$ are $\R$-polarizable.

  The {\it Mumford--Tate group $G$ of $H$} is 
the Tannaka group (cf.\ \cite{Mi}) of the Tannaka category generated by $H$ 
(cf.\ \cite{A}).
  Explicitly, it is the smallest $\Q$-subgroup $G$ of $\Aut(H_\Q)$ such that $G_\R$ contains the image of the homomorphism $h:S_{\C/\R} \to \Aut(H_\R)$ and such that $\Lie(G)_\R$ contains $\delta$. 
  Here $h$ and $\delta$ are determined by the canonical splitting of $H$ (\cite{CKS}, \cite{KNU} Part II 1.2).
  In the case where $H$ is pure, 
$G$ is the smallest $\Q$-subgroup of $\Aut(H_\Q)$ such that $G_\R$ contains the image of $S_{\C/\R}\to \Aut(H_\R)$. 

  In the pure case, the following proposition is well-known. 

\end{sbpara}

\begin{sbprop}\label{equiv} 
The two definitions of G in 4.2.1 coincide.
\end{sbprop}

\begin{pf} 
Let $G$ be the group explicitly defined in the latter part of 
\ref{MTgp}.

Let $J$ be the Tannaka group defined in the former part of \ref{MTgp}. 
By Tannaka duality, the theory of $(h, \delta)$ of MHS gives  a homomorphism $S_{\C/\R}\to J_\R$ and $\delta\in \Lie(J)_\R$. The homomorphism $J\to \Aut(H_\Q)$ is injective (otherwise, if $K$ denotes the kernel, representations of $J/K$ would form a smaller Tannaka category). 
Hence we have $G\subset J$. We will use

{\bf Claim.} {\it For linear representations $V_1$ and $V_2$ of $J$ over $\Q$, we have $\Hom_J(V_1, V_2)= \Hom_G(V_1, V_2)$. }

This is because $$\Hom_J(V_1, V_2)\subset \Hom_G(V_1,V_2)\subset \Hom_{\text{MHS}}(V_1, V_2)= \Hom_J(V_1, V_2).$$

By the pure case, we have $G/G_u = J/J_u$. (For this, a point is that $G_u$ coincides with the kernel $G_1$ of $G\to \Aut(\gr^W)$. We have $G_1\subset G_u$. It is sufficient to prove that $G/G_1$ is reductive. This is seen from the polarizability of $\gr^W$.)  

Assume $G\neq J$. Then by $G/G_u=J/J_u$, we have $G_u \neq J_u$. Hence the map $G_u\to J_u/[J_u, J_u]$ is not surjective. Since the image of this map is stable under the adjoint action of $G/G_u=J/J_u$, the image is a normal subgroup of $J/[J_u, J_u]$. Let $Q$ be the quotient of $J/[J_u, J_u]$ by this image. 
Let $Q_1$ be the quotient of $J_u/[J_u, J_u]$ by the image of $G_u$. Then $Q$ is a semi-direct product of $Q_1$ and $G/G_u$. We consider the following representations $V_1$ and $V_2$ of $Q$ over $\Q$. Let $V_1=\Q$ with the trivial action of $Q$. Let $V_2=\Q\oplus Q_1$ on which $G/G_u$ acts by the trivial action on $\Q$ and by the adjoint action on $Q_1$, and $v\in Q_1$ acts by sending $(1,0)$ to $(1, v)$ and trivially on $Q_1$. The $\Q$-linear map $V_1\to V_2$ which sends $1$ to $(1,0)$ is a $G$-homomorphism but not a $J$-homomorphism. This contradicts the Claim. 
\end{pf}

\begin{sbpara}\label{MTdom}

  The {\it Mumford--Tate domain associated to $H$} is defined as the period domain $D$ associated to $G$ and $h_0: S_{\C/\R}\to (G/G_u)_\R$ 
which is defined by $\gr^WH$.

  In the pure case, our $\Gamma \bs D_{\Sig}$ is essentially the same as the one by Kerr--Pearlstein (\cite{KP}).

\end{sbpara}

\subsection{Mixed Shimura varieties}
  See \cite{Mi} for the generality of mixed Shimura varieties. 
  This is the case where the universal object satisfies Griffiths transversality. $\gr^W_w\Lie(G)$ should be $0$ unless $w=0, -1, -2$. The $(p, q)$-Hodge component of $\gr^W_w\Lie(G)$ for $w=0$ (resp.\ $w=-1$, resp.\ $w=-2$) should be $0$ unless $(p,q)$ is $(1,-1)$, $(0,0)$, and $(-1,1)$ (resp.\ $(0, -1)$ and $(-1, 0)$, resp.\ $(-1,-1)$). (If this condition is satisfied by one point of $D$, it is satisfied by all points of $D$.)

For example, the universal abelian variety over a Shimura variety of 
PEL (polarizations, endomorphisms, and level structures) type is a mixed Shimura variety. Toroidal compactifications of 
these universal abelian varieties are expressed as $\Gamma \bs 
D_{\Sig}$.

\section{Higher Albanese manifolds and their toroidal partial compactifications}
\label{s:hab}

\subsection{Understanding higher Albanese manifolds by $D$}
\label{s:hab1}

\begin{sbpara}\label{hab1} 
Let $X$ be a connected smooth algebraic variety over $\C$. 
Fix $b\in X$. 
Let $\Gamma$ be a quotient group of $\pi_1(X, b)$ and assume that 
$\Gamma$ is a torsion-free nilpotent group.

Let $\cG=\cG_{\Gamma}$ be the unipotent algebraic group over $\Q$ 
whose Lie algebra is defined as follows. Let $I$ be the augmentation 
ideal $\text{Ker}(\bQ[\Gamma]\to \bQ)$ of $\Q[\Gamma]$. 
Then 
$\Lie(\cG)$ is the $\Q$-subspace of 
$\Q[\Gamma]^{\wedge}:=\varprojlim_n \Q[\Gamma]/I^n$ generated by all 
$\log(\gamma)$ ($\gamma \in \Gamma$). The Lie product of $\Lie(\cG)$ 
is defined by $[x,y]= xy-yx$. 
We have $\Gamma \subset \cG(\Q)$. 
 
We have 
$$\Lie(\cG)= \{h \in \Q[\Gamma]^{\wedge}\;|\; \Delta(h)= h \otimes 1 + 
1 \otimes h\},$$ 
$$\cG (R)= \{g\in (R[\Gamma]^{\wedge})^\times\;|\; \Delta(g)= g \otimes g\}$$ 
for any commutative ring $R$ over $\Q$, where $\Delta: 
R[\Gamma]^{\wedge}\to R[\Gamma\times \Gamma]^{\wedge}$ is the ring 
homomorphism induced by the ring homomorphism $R[\Gamma]\to 
R[\Gamma\times \Gamma]\;;\;\gamma \mapsto \gamma\otimes \gamma$ 
($\gamma \in \Gamma$).

\end{sbpara} 
 
\begin{sbpara}\label{hab-h}  
For $r\geq 0$, let $\Gamma_r$ be the 
torsion-free nilpotent quotient group of $\pi_1(X,b)$ defined in 
\ref{result1}. 
Then for a given $\Gamma$ as in \ref{hab1}, there is an 
$r\geq 1$ such that $\Gamma$ is a quotient of $\Gamma_r$. 
We define 
the weight filtration on $\Lie(\cG_{\Gamma})$ (resp.\ the Hodge 
filtration on $\Lie(\cG_{\Gamma})_\C$) as the image of that of 
$\Lie(\cG_{\Gamma_r})$ (resp.\ $\Lie(\cG_{\Gamma_r})_\C$). 
This gives a 
structure of an MHS on $\Lie(\cG_{\Gamma})$ which is independent of the 
choice of $r$. 
 
Note that $\cG_{\Gamma_r}$ is written as $\hat G_r$ in \ref{hamHZ}. 
\end{sbpara}

\begin{sbpara}\label{hab3} 
The higher Albanese manifold $A_{X,\Gamma}$ of $X$ for $\Gamma$ is as follows. 
Let $F^0\cG(\C)$ be the algebraic subgroup of $\cG(\C)$ over $\C$ corresponding to the Lie subalgebra $F^0\Lie(\cG)_\C$ of $\Lie(\cG)_\C$. Define
$$A_{X,\Gamma}:= \Gamma\bs \cG(\C)/F^0\cG(\C).$$

Let $\Gamma_r$ be as in \ref{result1}.
For $\Gamma=\Gamma_r$, $A_{X,\Gamma}$ coincides with $\Alb^r(X)$ in \ref{hamHZ}.

In the case where $\Gamma$ is $H_1(X, \Z)/(\text{torsion})$ regarded as a quotient group of $\pi_1(X,b)$, $A_{X,\Gamma}$ coincides with the Albanese variety 
$\Gamma\bs H_1(X,\C)/F^0H_1(X,\C)$ of $X$. 

We will give an understanding of $A_{X,\Gamma}$ by using $D$ of this paper in Theorem \ref{thm1}. 

We will describe the functor represented by $A_{X,\Gamma}$ in Theorem \ref{thm2}. 

\end{sbpara}

\begin{sbpara}\label{und2} Take a $\Q$-MHS $V_0$ with polarizable $\gr^W$ having the $\Q$-MHS $\Lie(\cG)$ as a direct summand.
Let $Q$ be the Mumford--Tate group (\ref{MTgp}) associated to the $V_0$
 (5.1.1).
The action of $Q$ on $\Lie(\cG)$ induces an action of $Q$ on $\cG$. 
By using this action, define the semidirect product $G$ of $Q$ and $\cG$ with an exact sequence $1\to \cG\to G\to Q\to 1$. We have $\cG\subset G_u$. We have $h_0: S_{\C/\R}\to (Q/Q_u)_\R= (G/G_u)_\R$ given by the Hodge decomposition of $\gr^WV_0$.

Then $(G, \Gamma)$ satisfies the condition in \ref{Gamma}, and $\Gamma$ is a neat subgroup of $G(\Q)$. 

Let $D_G$ (resp.\ $D_Q$) be the period domain $D$ for $G$ (resp.\ $Q$) and $h_0$ (\ref{pd}). 
We have a canonical map $\Gamma \bs D_G\to D_Q$ induced by the canonical homomorphism $G\to Q$. 

\end{sbpara}

\begin{sbpara}\label{und8}
We define $b_Q\in D_Q$ and $b_G\in D_G$ as follows.
Let $b_Q\in D_Q$ be the isomorphism class of the evident functor $\Rep(Q)\to\text{$\Q$-MHS}$ (cf.\ Proposition \ref{hdelta}).
Since $Q\subset G$ is a semidirect summand, we have the restriction functor of $\Rep(G)\to\Rep(Q)$.
Let $b_G\in D_G$ be the isomorphism class of the composite functor $\Rep(G)\to\Rep(Q)$ and $b_Q:\Rep(Q)\to\text{$\Q$-MHS}$.
Then we see that the map $D_G\to D_Q$, induced by the canonical homomorphism $G\to Q$, sends $b_G$ to $b_Q$.
Let $\cD$ be the fiber of the map $D_{G} \to D_{Q}$ over $b_Q$.

The following theorem is a generalization of \cite{HZ1} (5.10) into the present context of tensor functors.

\end{sbpara}

\begin{sbthm}\label{thm1} 
The map $\cG(\C) \to D_G\; ;\; g \mapsto gb_G$ induces isomorphisms$:$

$(1)$ $\cG(\C)/F^0\cG(\C) \cong \cD$. 

$(2)$ $A_{X,\Gamma}\cong \Gamma \bs \cD$. 
\end{sbthm}

\begin{pf} 
We prove (1), from which (2) follows.
Define
$$
F^0(G(\R)G_u(\C)):=\{g\in G(\R)G_u(\C)\,|\, gb_G=b_G\},
$$
$$
F^0(Q(\R)Q_u(\C)):=\{g\in Q(\R)Q_u(\C)\,|\, gb_Q=b_Q\}.
$$
Then we have a commutative diagram of exact sequences
$$
\CD
1@>>>\cG(\C)@>>>G(\R)G_u(\C)@>>>Q(\R)Q_u(\C)@>>>1\\
&&\bigcup&&\bigcup&&\bigcup\\
1@>>>F^0(\cG(\C))@>>>F^0(G(\R)G_u(\C))@>>>F^0(Q(\R)Q_u(\C))@>>>\,1.
\endCD
$$
Here the surjectivity of $F^0(G(\R)G_u(\C))\to F^0(Q(\R)Q_u(\C))$ follows from 
$F^0(Q(\R)Q_u(\C))\subset F^0(G(\R)G_u(\C))$ which is induced from $Q(\R)Q_u(\C)\subset G(\R)G_u(\C)$.

Combining this with $(G(\R)G_u(\C))/F^0\overset{\sim}{\to} D_G$; $g\mapsto gb_G$ and $(Q(\R)Q_u(\C))/F^0\overset{\sim}{\to} D_Q$; $g\mapsto gb_Q$, we get $ \cG(\C)/F^0\overset{\sim}{\to} \cD$; $g\mapsto gb_G$.
\end{pf}

\begin{sbpara}\label{hab2} 
Let $\cC_{X,\Gamma}$ be the category of variations of $\Q$-MHS $\cH$ 
on $X$ satisfying the following conditions. 
 
(i) For any $w\in \Z$, 
$\gr^W_w\cH$ is a constant  polarizable Hodge structure. 
 
(ii) $\cH$ is good at infinity in the sense of \cite{HZ1} (1.5). 
 
(iii) The monodromy action of $\pi_1(X, b)$ on $\cH_{\Q}(b)$ (which is 
unipotent under (i)) factors through $\Gamma$.

Let $\cC_{X,\Gamma}'$ be the category of $\Q$-MHS $H$ with polarizable 
$\gr^W$ endowed with an action of the Lie algebra $\Lie(\cG)$ on 
$H_{\Q}$ such that  $\Lie(\cG) \otimes H \to H$ is a homomorphism of 
MHS. 
 
\end{sbpara} 
 
\begin{sbpara}\label{HZgam} 
 Let $\cC_{X, r}$ (resp.\ $\cC'_{X,r}$) be the left (resp.\ right) 
category in the equivalence of categories at the beginning of 1.1. 
Then we have 
$$\cC_{X,\Gamma_r} \supset \cC_{X,r}, \quad \cC'_{X,\Gamma_r}\supset 
\cC'_{X,r},$$ 
$$\bigcup_{\Gamma} \; \cC_{X,\Gamma}=\bigcup_r \; \cC_{X, 
\Gamma_r}=  \bigcup_r \; \cC_{X,r}, \quad \bigcup_{\Gamma}\;  \cC'_{X, 
\Gamma}=\bigcup_r\;  \cC'_{X, \Gamma_r}=\bigcup_r\;  \cC'_{X,r}.$$

The equivalences $\cC_{X,r}\simeq \cC'_{X,r}$  in 1.1 for all $r$ 
induce an equivalence $\bigcup_r \; \cC_{X,r}\simeq \bigcup_r \; 
\cC'_{X,r}$, and it induces an equivalence  $$\cC_{X,\Gamma}\simeq 
\cC'_{X,\Gamma}$$ between the full subcategories. 
 
\end{sbpara}

\begin{sbpara}\label{moduli1} 
Define a contravariant functor 
$$
\cF_{\Gamma}:\cB(\log) \to(\text{Set})
$$
as follows. 

${\cF}_{\Gamma}(S)$ is the set of isomorphism classes of pairs $(H,\mu)$, where $H$ is an exact $\otimes$-functor
$\cC_{X,\Gamma} \to \text{MHS}(S)$ and  $\mu$ is a $\Gamma$-level structure, satisfying the following 
condition (i). 
Here a {\it $\Gamma$-level structure} means a global section 
of the sheaf $\Gamma \bs \Cal I$, where $\Cal I$ is the sheaf of functorial 
$\otimes$-isomorphisms $H(\cH)_{\Q} \overset{\cong}\to \cH(b)_{\Q}$ of $\Q$-local systems preserving weight filtrations.

 (i) For any $\Q$-MHS $h$, we have a functorial $\otimes$-isomorphism $H(h_X) \cong h_S$ such that the induced isomorphism of local systems $H(h_X)_\Q\cong h_\Q=h_X(b)_\Q$ belongs to $\mu$. Here $h_X$ (resp.\  $h_S$) denotes the constant variation (resp.\ family) of $\Q$-MHS over $X$ (resp.\ $S$) associated to $h$. 
 
\end{sbpara}

\begin{sbthm}\label{thm2} 
The higher Albanese manifold $A_{X,\Gamma}$ represents  $\cF_{\Gamma}$. 
\end{sbthm}

\begin{pf}
For $S\in \cB(\log)$, we show $\Mor(S,A_{X,\Gamma})\simeq\cF_\Gamma(S)$.

The map from the right-hand-side to the left-hand-side is as follows.
For an element $\cC_{X,\Gamma}\to\text{MHS}(S)$ of $\cF_{\Gamma}(S)$, consider the composition
$$
\Rep(G)\subset\cC'_{X,\Gamma}\simeq\cC_{X,\Gamma}\to\text{MHS}(S).
$$
By the non-log version of the general theorem \ref{t:extend}, this yields a morphism $S\to\Gamma\bs D_G$ whose image is sent to $b_Q$ under $\Gamma\bs D_G\to  D_Q$.
Thus we get an element $S\to\Gamma\bs\cD=A_{X,\Gamma}$ (\ref{thm1}) of $\Mor(S,A_{X,\Gamma})$.

As for the map from the left-hand-side to the right-hand-side, which is inverse to the above map, we give two constructions.

The first construction is as follows.
Assume that we are given a morphism $S \to A_{X, \Gamma}$. 
Similarly as in \ref{inv}, for an object $V$ of $\cC'_{X,\Gamma}$, we have a Lie algebra homomorphism $\Lie(\cG)\to\End(V)$ which is a homomorphism of MHS, and it induces a morphism from $A_{X,\Gamma}$ to the classifying space $\Gamma\bs D(V)$ for $V$.
Pulling back the universal variation of MHS on $\Gamma\bs D(V)$ by the composition 
$S\to A_{X,\Gamma}\to\Gamma\bs D(V)$, we get an object of $\text{MHS}(S)$.
This gives a desired pair $(H,\mu)$ of a functor $H:\cC_{X,\Gamma}\simeq\cC'_{X,\Gamma}\to\text{MHS}(S)$ and a $\Gamma$-level structure $\mu$.

The second construction is as follows.
Assume we are given a morphism $S\to A_{X,\Gamma}$ and an object $V$ of $\cC'_{X,\Gamma}$. 
Let $Q$ be the Mumford--Tate group of the $\Q$-MHS $V_0:=\Lie(\cG)\oplus V$, and define $G$ as in \ref{und2}. 
Then we have
$S\to A_{X,\Gamma}\simeq\Gamma\bs\cD\hookrightarrow \Gamma\bs D_G$.
By the non-log version of Theorem \ref{t:main}, the object $V$ of $\text{Rep}(G)$ gives an 
object of MHS(S).\end{pf}

\begin{sbpara}

The {\it higher Albanese map} $\varphi:X\to A_{X,\Gamma}$ 
corresponds in \ref{thm2} to  the evident functor $H:\cC_{X,\Gamma} 
\to \text{MHS}(X)$. 
 
\end{sbpara}

\subsection{Toroidal partial compactifications}\label{s:tpc}

\begin{sbpara}\label{tpc}
Let $G$ be as in \ref{und2}.
Let $\Sig$ be a weak fan in $\Lie(G)$ such that $\sig\subset \Lie(\cG)_\R$ for any $\sig\in \Sig$. 
Assume that $\Sigma$ and $\Gamma$ in \ref{hab1} are strongly compatible (\ref{fan}).
Then, we have a canonical morphism $\Gamma \bs D_{G,\Sig}\to D_Q$, extending the morphism $\Gamma \bs D_G\to D_Q$, induced by the homomorphism $G\to Q$.

\end{sbpara}

\begin{sbpara}
Define the toroidal partial compactification $
{A}_{X, \Gamma,\Sig}$ of $A_{X,\Gamma}$ as the subspace of $\Gamma \bs D_{G,\Sig}$ which is defined to be the inverse image of $b_Q$. We can endow $
{A}_{X, \Gamma, \Sig}$ with a structure of a log manifold such that for any object $S$ of $\cB(\log)$,  $\text{Mor}(S, 
{A}_{X,\Gamma,\Sig})$ coincides with the set of all morphisms $S \to \Gamma \bs D_{G,\Sig}$ whose images in $D_Q$ are $b_Q$. 

This $A_{X, \Gamma, \Sig}$ is independent of the choice of $V_0$ in \ref{und2} which is used in the definitions of $Q$ and $G$.

\end{sbpara}

\begin{sbpara}\label{moduli} 
Define a contravariant functor 
$$\cF_{\Gamma,\Sig}:\cB(\log) \to(\text{Set})$$
as follows.

${\cF}_{\Gamma,\Sig}(S)$ is the set of isomorphism classes of pairs $(H,\mu)$ where $H$ is an 
exact $\otimes$-functor $\cC_{X,\Gamma}\to\LMH(S)$
and  $\mu$ is a $\Gamma$-level structure
satisfying the condition (i) in \ref{moduli1}
and also the following condition (ii).
 
 (ii) The following (ii-1) and (ii-2) are satisfied for any $s\in S$ and any $t\in s^{\log}$. Let $\tilde \mu_t: H(\cH)_{\Q,t}\cong \cH(b)_\Q$ be a  functorial $\otimes$-isomorphism which belongs to $\mu_t$.

  (ii-1) There is a $\sig\in \Sig$ such that the logarithm of the action of the local monodromy cone 
$\Hom((M_S/\cO^\times_S)_s, \N)\subset \pi_1(s^{\log})$ on $H_{\Q,t}$ is contained, via $\tilde \mu_t$, in $\sigma \subset \Lie(\cG)_\R$.

  (ii-2) Let $\sig\in \Sig$ be the smallest cone which satisfies (ii-1) and let $a: \cO_{S,t}^{\log}\to \C$ be a ring homomorphism which induces the evaluation $\cO_{S,s}\to \C$ at $s$. Then, for each $\cH\in \cC_{X,\Gamma}$, $(\sig, \tilde \mu_t(a(H(\cH))))$ generates a nilpotent orbit in the sense of \cite{KNU} Part III, 2.2.2.   
  
 \end{sbpara}

\begin{sbthm}\label{thm3} 
The functor $
{\cF}_{\Gamma,\Sig}$ is represented by $
{A}_{X, \Gamma,\Sig}$.
\end{sbthm}

\begin{pf}
This follows from Theorem \ref{t:main}.
The proof is similar to the one of Theorem \ref{thm2}. We replace $\text{MHS}(S)$ there by $\text{LMH}(S)$ and, in the latter half of the proof, we use the second construction of the inverse map. 
\end{pf}

\begin{sbpara} 
Let $\Xi$ be the set of all rational nilpotent cones in $\Lie(\cG)_{\R}$ of 
rank $\leq 1$. 
Then $\Xi$ is a fan and is strongly compatible with 
$\Gamma$. 
 
\end{sbpara}

\begin{sbthm}\label{a:extend} 
Let $\overline{X}$ be a smooth algebraic variety over $\C$ which contains $X$ 
as a dense open subset such that the complement 
$\overline{X}\smallsetminus X$ is a smooth divisor. 
Endow $\overline{X}$ with the log structure associated to this divisor. 
 
Then the higher Albanese map $\varphi: X \to A_{X,\Gamma}$ extends 
uniquely to a morphism $\overline{\varphi}: \overline{X} \to 
A_{X,\Gamma,\Xi}$ of log manifolds giving a commutative diagram 
$$ 
\CD 
X&\subset&\overline{X}\\ 
@V{\varphi}VV@VV{\overline{\varphi}}V\\ 
A_{X,\Gamma}&\quad \subset\quad& \, A_{X,\Gamma,\Xi}. 
\endCD 
$$ 
\end{sbthm}

\begin{pf} 
Since an object of $\cC_{X, \Gamma}$ is good at infinity, it extends to an LMH over $\overline{X}$.
Hence this theorem follows from (2) of the the general theorem \ref{t:extend} by using 
Theorem  \ref{thm3}.
\end{pf}

\subsection{Example}\label{P1-3}

\begin{sbpara}\label{E1}
We consider $$X:= {\bf P}^1(\C)\smallsetminus \{0,1,\infty\}\subset \overline{X}:= {\bf P}^1(\C).$$ We will consider the  toroidal partial compactification of the second higher Albanese manifold $\Alb^2(X)$ (\ref{hamHZ}) and the extended higher Albanese map from $\overline{X}$ to it (\ref{a:extend}).  The description of the degeneration  at the boundary of $\overline{X}$ becomes simpler if we take the base point 
$b$ of the theory of Hain--Zucker in the boundary outside $X$. 
For this, we can use the idea of tangential base point of Deligne (\cite{De3} Section 15 \lq\lq Points base \`a l'infini") and  its variant described in \ref{bbase}. 

As in \cite{HZ2}, 
$$\Alb^2(X) \cong \begin{pmatrix}  1 & \Z & \Z\\ 0 & 1 & \Z\\ 0 & 0 & 1 \end{pmatrix}\left\backslash \begin{pmatrix}  1 & \C& \C\\ 0 & 1 & \C\\ 0 & 0 & 1 \end{pmatrix}\right.$$
The right-hand-side is actually a period domain $G_{u,\Z}\bs D(\La)$ of classical type (see \ref{E11} below) and the toroidal partial compactification of $\Alb^2(X)$ with respect to the fan $\Xi$ (5.2.5) is isomorphic to the toroidal partial compactification $G_{u,\Z}\bs D(\La)_{\Xi}$ of this classical period domain considered in \cite{KNU} Part III. We first consider this period domain of the classical type in \ref{E11}--\ref{E14}. 

\end{sbpara}

\begin{sbpara}\label{E11}
Let $\Lam= (H_0, W, (\langle \;,\;\rangle)_{w\in \Z}, (h^{p,q})_{p,q\in \Z})$ be as follows. 
$H_0$ is a free $\Z$-module of rank $3$ with basis $(e_j)_{1\leq j\leq 3}$, $W$ is the increasing filtration on $H_{0,\Q}$ defined as
$$W_{-5}=0\subset W_{-4}=W_{-3}= \Q e_1\subset W_{-2}=W_{-1}=\Q e_1+\Q e_2\subset W_0=H_{0,\Q},$$
$\langle\;,\;\rangle_w : \gr^W_w(H_{0,\Q}) \times \gr^W_w(H_{0,\Q})\to \Q$  are the $\Q$-bilinear forms characterized by $\langle e_3,e_3\rangle_0=\langle e_2,e_2\rangle_{-2}= \langle e_1,e_1\rangle_{-4}=1$, and $h^{0,0}=h^{-1,-1}=h^{-2,-2}=1$, $h^{p,q}=0$ for the other $(p,q)$. 
For $R=\Z, \C$, let $G_{u,R}$ be the group of  automorphisms of the $R$-module $H_{0,R}$ which preserve $W$ and induce the identity map on $\gr^W$. Then $G_{u,R}$ is identified with the group of unipotent upper triangular $(3,3)$-matrices with entries in $R$. 

The period domain $D(\Lam)$ is isomorphic to $G_{u,\C}$ where
the matrix 
$$\begin{pmatrix}  1 & \beta & \lambda\\ 0 & 1 & \alpha \\ 0 & 0 & 1\end{pmatrix}$$
corresponds to the following decreasing filtration $F=F(\alpha,\beta, \lam)$ on $H_{0,\C}$: 
$F^1=0$, $F^0$ is generated by $e_3+\alpha e_2+ \lam e_1$, $F^{-1}$ is generated by $F^0$ and $e_2+\beta e_1$, and $F^{-2}=H_{0,\C}$.

The natural  action of $G_{u,\C}$ on $D(\La)$ is identified with the natural action of $G_{u,\C}$ on itself from the left. 

\end{sbpara}

\begin{sbpara}\label{E12} Let $\Xi$ be the set of all cones of the form $\R_{\geq 0} N$ where $N$ is a $\Q$-linear map $H_{0,\Q}\to H_{0,\Q}$ such that $NW_w\subset W_{w-1}$ for all $w\in \Z$. 
We consider the extended period domain 
$D(\Lam)_{\Xi}$ (\cite{KNU} Part III).

For $N: H_{0,\Q}\to H_{0, \Q}$ defined by $Ne_3= ae_2+ce_1$, $Ne_2=be_1$, $Ne_1=0$ ($a,b, c\in \Q$), $(N, F(\alpha, \beta, \lam))$ generates a nilpotent orbit if and only if it satisfies the Griffiths transversality $NF^p\subset F^{p-1}$ ($p\in \Z$), and hence if and only if $c= a\beta- b\alpha$. Hence for $\sig=\R_{\geq 0}N$, 
$D_{\sig}\neq D$ if and only if either $a\neq 0$ or $b\neq 0$. 
\end{sbpara}

\begin{sbpara}\label{E13}

By \cite{KNU} Part III,  
the quotient $G_{u,\Z}\bs D(\Lam)_{\Xi}$ has a structure of a log manifold. 

Let $N:H_{0,\Q}\to H_{0,\Q}$ be the $\Q$-linear map defined by $Ne_3=e_2$, $Ne_2=Ne_1=0$ and let $\sig=\R_{\geq 0} N$. 

We describe the local structure of $ G_{u,\Z}\bs D(\Lam)_{\Xi}$ at the point corresponding to a $\sig$-nilpotent orbit. Let $p$ be the image of an element of $D(\La)_{\sig}\smallsetminus D(\La)$ in $G_{u,\Z}\bs D(\La)_{\Xi}$. Then for some $\lam_0\in \C$, $p$ is the class of the $\sig$-nilpotent orbit generated by $(N, F(0,0,\lambda_0))$. Let $Y$ be the log manifold $\{(q, \beta, \lam)\in \C^3\;|\; \beta=0 \;\text{if}\; q=0\}$ with the strong topology (\cite{KU} Section 3.1), with the structure sheaf of rings which is the inverse image of the sheaf of holomorphic functions on $\C^3$, and with the log structure generated by $q$. Then there is an open neighborhood $U$ of $(0,0,\lam_0)$ in $\C^3$ and an open immersion $$Y\cap U \overset{\subset}\to G_{u,\Z}\bs D(\La)_{\Xi}$$ of log manifolds which sends $(q, \beta, \lam)\in Y\cap U$ with $q\neq 0$ to the class 
of $F(\alpha, \beta, \lam)=\exp(\alpha N)F(0, \beta, \lam)$, where $\alpha\in \C$ is such that $q=e^{2\pi i\alpha}$, and which sends $(0,0,\lam_0)$  to $p$. 
\end{sbpara}

\begin{sbpara}\label{E14} We can show that for any $p\in G_{u,\Z}\bs D(\La)_{\Xi}$ which does not belong to $G_{u,\Z}\bs D(\La)$, there are an open neighborhood $U$ of $(0,0,0)$ in $\C^3$ and an open immersion $Y\cap U \to G_{u,\Z}\bs D(\La)_{\Xi}$ of log manifolds which sends $(0,0,0)$ to $p$. 

\end{sbpara}

\begin{sbpara}\label{bbase} We describe how to formulate  a base point in the boundary in the theory of Hain--Zucker. The following is a variant of tangential base point of Deligne and matches log Hodge theory well.

Let $\overline{X}$ be a connected smooth algebraic variety over $\C$, let $D$ be a divisor on $\overline{X}$ with normal crossings, and let $X:=\overline{X}\smallsetminus D$. Endow $\overline{X}$ with the log structure associated to $D$. We formulate  a base point in $\overline{X}$ outside $X$ as follows. 

In our definition, a base point in the boundary of $\overline{X}$ is a pair $b=(y, a)$ where $y$ is a point of $\overline{X}^{\log}$ which does not belong to $X$, and $a$ is a specialization $\cO^{\log}_{\overline{X},y}\to \C$. That is, $y$ is a pair $(x,h)$ where $x$ is a point of $\overline{X}$ which does not belong to $X$, $h$ is a homomorphism $M^{\gp}_{\overline{X},x} \to {\bf S}^1:=\{z\in \C^\times \;|\; |z|=1\}$ whose restriction to the subgroup $\cO_{\overline{X},x}^\times$ of $M_{\overline{X},x}^{\gp}$ coincides with $f \mapsto f(x)/|f(x)|$, and $a$ is 
a ring homomorphism $\cO^{\log}_{\overline{X},y}\to \C$ whose restriction to the  subring $\cO_{\overline{X},x}$ of $\cO^{\log}_{\overline{X}, y}$ coincides with $f\mapsto f(x)$. 

A path between $y$ and a point $b'$ of $X$ induces an isomorphism $$\pi_1(\overline{X}^{\log}, y) \cong \pi_1(X,b').$$ We can use $\pi_1(\overline{X}^{\log}, y)$ in place of $\pi_1(X,b')$ ($b'\in X$) in the theory of Hain--Zucker.

Let $\Gamma$ be a nilpotent torsion free quotient group of $\pi_1(\overline{X}^{\log}, y)$. Then we have the unipotent algebraic group $\cG$ over $\Q$ by using $\pi_1(\overline{X}^{\log}, y)$ by the method of \ref{hab1}. We obtain a structure of a $\Q$-MHS on $\Lie(\cG)$ as follows. For $b'\in X$, let $\cG(b')$ be the unipotent group $\cG$ in \ref{hab1} obtained by using the base point $b'$. Then when $b'\in X$ moves, the MHS $\Lie(\cG(b'))$ forms an object $\cH$ of $\cC_{X,\Gamma}$. Since $\cH$ is good at infinity, it extends uniquely 
to a $\Q$-LMH $\overline{\cH}$ on $\overline{X}$. By the  specialization of $\overline{\cH}$ by $a$ at $y$, we obtain a structure of $\Q$-MHS on $\Lie(\cG)$. 

The main theorem of Hain--Zucker \cite{HZ1} introduced in Section 1.1 and results in Sections 5.1 and 5.2 remain true when we use the base point in the boundary,  and are  deduced from the work \cite{HZ1} and by the arguments in Sections 5.1 and 5.2. 
\end{sbpara}

\begin{sbpara}\label{tbase}

A tangential base point of Deligne 
 in the boundary of $\overline{X}$ (\cite{De3} Section 15) gives a base point $b$ in the boundary in our sense (\ref{bbase}). We explain this in the case where $X$ is a curve. In this case, a tangential base point in the boundary is a non-zero element $v$ of the tangent space $T_x(\overline{X})=\Hom_\C(m_x/m_x^2, \C)$ with $x\in \overline{X}\smallsetminus X$ and $m_x$ being the maxmal ideal of $\cO_{\overline{X},x}$. 
We have the  corresponding base point $b=(y, a)$, $y=(x,h)$ in the boundary in our sense  as follows.

$h:M^{\gp}_x \to {\bf S}^1$ is the unique homomorphism which sends any element $f$ of $\cO_{\overline{X},x}^\times$ to $f(x)/|f(x)|$ and any prime element $t$ of $\cO_{\overline{X},x}$ to $v(t)/|v(t)|$. 
$a$ is the unique ring homomorphism $\cO^{\log}_{\overline{X}, y}\to \C$ satisfying the following (i) and  (ii).

(i) $a(f)=f(x)$ for any $f\in \cO_{\overline{X},x}$.

(ii) Let $t$ be a prime element of $\cO_{\overline{X}, x}$ such that $h(t)=1$ (that is, $v(t) \in \R_{>0}$). Let $f\in \cO^{\log}_{\overline{X}, y}$ be the branch of $\log(t)$ such that the imaginary part of $f(x')$ converges to $0$ if $x'\in X$ converges to $y$ in $\overline{X}^{\log}$. Then $a(f)= \log(v(t))\in \R$.

\end{sbpara}

\begin{sbpara}\label{E21}
Now let 
$X= {\bf P}^1(\C)\smallsetminus \{0,1,\infty\}$, $\overline{X}={\bf P}^1(\C)$. We take  the base point $b$ in the boundary of $\overline{X}$ corresponding to the tangent vector $v$ at $0\in \overline{X}$ which sends the class of the coordinate function $x$ of $\C\subset {\bf P}^1(\C)$ in $m_0/m_0^2$ to $1$. That is, 
$b= (y, a)$, $y=(0,h)\in \overline{X}^{\log}$ where $h$ sends the coordinate function $x$ to $1$ and $a$ sends the branch of $\log(x)$ which has real value on $\R_{>0}$ to $0$.

The group $\pi_1(\overline{X}^{\log}, y)$ is a free group of rank $2$ generated by elements $\gamma_0$ and $\gamma_1$. Here for $\alpha=0,1$, $\gamma_{\alpha}$ is the class of the following  loop $[0,1]\to \overline{X}^{\log}$ which we denote also by $\gamma_{\alpha}$.
Let $x$ be the coordinate function of $\C\subset {\bf P}^1(\C)$. 
Then $\gamma_0(t)=(0, h)$ where $h$ sends $x$ to $e^{2\pi it}$. 
$\gamma_1(0)=\gamma_1(1)=(0,h)$ where $h$ sends $x$ to $1$.
For $0<t<1/3$, $\gamma_1(t)= 3t\in X$. For $1/3\leq t\leq 2/3$, $\gamma_1(t)= (1, h)$ where $h$ sends $1-x$ to $e^{2\pi i (3t-1)}$. For $2/3<t<1$, $\gamma_1(t)= 3(1-t) \in X$. 
\end{sbpara}

\begin{sbpara}\label{E22}
Let $\Gamma$ be the quotient group $\pi_1(\overline{X}^{\log},y)/[\pi_1, [\pi_1, \pi_1]]$ of $\pi_1=\pi_1(\overline{X}^{\log}, y)$. 
 We consider $A_{X,\Gamma}$ (\ref{hab3}) which is the second higher Albanese manifold of $X$ by using the above base point $b$ in the boundary. 

For $\alpha=0, 1$, let $N_{\alpha}= \log(\gamma_{\alpha})\in \Lie(\cG)$. Then $\Lie(\cG)$ is three-dimensional over $\Q$ with basis $N_0$, $N_1$, $[N_1,N_0]$ (cf.\ \ref{hab1}).

The mixed Hodge structure on $\Lie(\cG)$ is as follows. The weight filtration is given by 
$$W_{-5}=0\subset W_{-4}=W_{-3}= \Q \cdot[N_1,N_0]\subset W_{-2}= \Lie(\cG).$$ $N_0$ and  $N_1$ are  of Hodge type $(-1,-1)$, and  $[N_1, N_0]$ is of Hodge type $(-2,-2)$.

We have $F^0\cG(\C)=\{1\}$ and hence 
$$A_{X,\Gamma}= \Gamma\bs \cG(\C).$$
\end{sbpara}

The following \ref{E23} and \ref{E24} are seen from \cite{De3} (cf. also \cite{HZ2}). 
\begin{sbpara}\label{E23}

Consider the following $\Q$-MHS $V$ and the Lie action of $\Lie(\cG)$ on $V$.
$V=H_{0,\Q}$ with the Hodge filtration $F(0,0,0)$ on $V_\C$ (\ref{E11}). 
The action of $\Lie(\cG)$ is as follows. $N_0e_3= e_2$, $Ne_j=0$ for $j=1,2$; $N_1e_2=e_1$, $N_1e_j=0$ for $j=1,3$. 
Then the action  $\Lie(\cG) \otimes V\to V$ is a homomorphism of MHS. 

This induces an isomorphism $A_{X,\Gamma}\overset{\cong}\to  G_{u,\Z}\bs D(\La)$ of complex analytic manifolds. It extends to an isomorphism $A_{X, \Gamma, \Xi} \cong G_{u,\Z}\bs D(\La)_{\Xi}$ of log manifolds. 

The composition $X\to G_{u,\Z}\bs D(\La)$ of the higher Albanese map $X\to A_{X,\Gamma}$ and the above isomorphism sends $x\in X$ to the class of $F((2\pi i)^{-1}\log(x), (2\pi i)^{-1}l_1(x), (2\pi i)^{-2}l_2(x))$, where $l_1(x)= -\log(1-x)$ and $l_2(x)$ is the dilog function.

\end{sbpara}

\begin{sbpara}\label{E24} 
There is another isomorphism $A_{X,\Gamma}\cong G_{u,\Z}\bs D(\La)$ which may be more popular. Consider the $\Q$-MHS $V$ as in \ref{E23} and consider the Lie action of $\Lie(\cG)$ on $V$ such that $N_0e_2=e_1$, $N_0e_j=0$ for $j=1,3$; $N_1e_3=e_2$, $N_1e_j=0$ for $j=1,2$.  Then the action  $\Lie(\cG) \otimes V\to V$ is a homomorphism of MHS. This induces an isomorphism $A_{X,\Gamma}\overset{\cong}\to  G_{u,\Z}\bs D(\La)$ of complex analytic manifolds and  an isomorphism $A_{X, \Gamma, \Xi} \cong G_{u,\Z}\bs D(\La)_{\Xi}$ of log manifolds.

In this case, the composition $X\to G_{u,\Z}\bs D(\La)$ of the higher Albanese map $X\to A_{X,\Gamma}$ and the above isomorphism sends $x\in X$ to the class of
$$\begin{pmatrix} 1 & -(2\pi i)^{-1}\log(x) & (2\pi i)^{-2}l_2(x)\\0 & 1 & -(2\pi i)^{-1}l_1(x)\\ 0&0&1\end{pmatrix}^{-1}\in G_{u,\C}\cong D(\La).$$

The pullback on $X$ of the universal object on $G_{u,\Z}\bs D(\La)$ under this composite map is the so-called dilog sheaf on $X$.

\end{sbpara}

\begin{sbpara}\label{E25} Consider the extended higher Albanese map $\overline{X} \to A_{X,\Gamma, \Xi}$ (\ref{a:extend}). Let $\overline{X} \to G_{u,\Z}\bs D(\La)_{\Xi}$ be the composite of this extended map and the isomorphism $A_{X,\Gamma,\Xi} \cong G_{u,\Z}\bs D(\La)_{\Xi}$ in \ref{E23}. Then the image of $0\in \overline{X}$ under this composite map is the class of the nilpotent orbit generated by $(N, F(0,0,0))$ with $N$ as in \ref{E13}. Let $Y\cap U \to G_{u,\Z} \bs D(\La)_{\Xi}$ be the open immersion given in \ref{E13} with $\lam_0=0$. Then if $x\in X$ is near to $0\in \overline{X}$, the image of $x$ in $G_{u,\Z}\bs D(\La)_{\Xi}$ is the image of $$(x, (2\pi i)^{-1}l_1(x), (2\pi i)^{-2}l_2(x))\in Y\cap U.$$ The last element converges to $(0,0,0)$ when $x$ converges to $0\in \overline{X}$. 

\end{sbpara}

{\bf Acknowledgments.} 
The authors thank Teruhisa Koshikawa for helpful discussions. 
The book \cite{Ko} of Toshitake Kohno was helpful for the authors. 
The authors thank to the referee for valuable comments. 
K.\ Kato was 
partially supported by NFS grants DMS 1001729 and DMS 1303421.
C.\ Nakayama was 
partially supported by JSPS Grants-in-Aid for Scientific Research (C) 22540011, (B) 23340008, and (C) 16K05093.
S.\ Usui was 
partially supported by JSPS Grants-in-Aid for Scientific Research (B) 23340008.

\noindent {\rm Kazuya KATO
\\
Department of mathematics
\\
University of Chicago
\\
Chicago, Illinois, 60637, USA}
\\
{\tt kkato@math.uchicago.edu}

\bigskip

\noindent {\rm Chikara NAKAYAMA
\\
Department of Economics 
\\
Hitotsubashi University 
\\
2-1 Naka, Kunitachi, Tokyo 186-8601, Japan
\\
{\tt c.nakayama@r.hit-u.ac.jp}

\bigskip

\noindent
{\rm Sampei USUI
\\
Graduate School of Science
\\
Osaka University
\\
Toyonaka, Osaka, 560-0043, Japan}
\\
{\tt usui@math.sci.osaka-u.ac.jp}


\begin{thebibliography}{99}
\bibitem{A} Y.\ Andr\'e,
{\rm Mumford-Tate groups of mixed Hodge structures and the theorem of the fixed part}, 
Comp.\ Math.\ {\bf 82-1} (1992), 1--24.

\bibitem{CKS} E.\ Cattani, A.\ Kaplan and W.\ Schmid,
{\rm Degeneration of Hodge structures},
Ann. of Math. {\bf 123} (1986), 457--535.

\bibitem{Ch} K.\ T.\ Chen, 
{\rm Extension of $C^{\infty}$ function algebra by integrals and Malcev compeletion of $\pi_1$}, 
Bull.\ of Amer.\ Math.\ Soc.\ {\bf 83} (1977), 831--879. 

\bibitem{De} P.\ Deligne, 
{\rm  Travaux de Shimura}, 
S\'eminaire Bourbaki (1970/71), Exp. 389, Lecture Notes in Math.\ {\bf 244}, Springer, (1971), 123--165.

\bibitem{De2} P.\ Deligne, 
{\rm La conjecture de Weil pour les surfaces K3}, 
Invent. Math. {\bf 15} (1972), 206--226. 

\bibitem{De3} P.\ Deligne, 
{\rm Groupe fondamental de la droite projective moins trois points, in Galois group over $\Q$}, 
MSRI Publications {\bf 16}  (1989), 79--297. 

\bibitem{DM} P.\ Deligne and J.\ S.\ Milne,  
{\rm Tannakian categories, in Hodge cycles, Motives, and Shimura varieties}, 
Lecture Notes in Math.\ {\bf 900}, Springer, (1982), 101--228. 


\bibitem{GGK} M.\ Green, P.\ Griffiths and M.\ Kerr, 
{\rm  Mumford-Tate Groups and Domains: Their Geometry and Arithmetic},
Ann.\ Math.\ Studies {\bf 183}, Princeton University Press (2012).

\bibitem{Gr} P.\ Griffiths,
{\rm  Periods of integrals on algebraic manifolds. I. Construction and properties of modular varieties}, 
Amer.\ J.\ Math.\ {\bf 90} (1968), 568--626.

\bibitem{He} R.\ M.\ Hain, 
{\rm The geometry of mixed Hodge structures on the fundamental group}, 
Proc.\ Symp.\ Pure Math.\ {\bf 46} (1987), 247--282. 

\bibitem{HZ1} R.\ M.\ Hain and S.\ Zucker,
{\rm Unipotent variation of mixed Hodge structure},
Invent.\ Math.\ {\bf 88} (1987), 83--124. 

\bibitem{HZ2} R.\ M.\ Hain and S.\ Zucker,
{\rm A guide to unipotent variations of mixed Hodge structures}, 
Lecture Notes in Math.\ {\bf 1246}, Springer, (2006), 92--106. 



\bibitem{KK} K.\ Kato and T.\ Koshikawa, 
Arithmetic models of Griffiths period domains, in preparation. 


\bibitem{KNU}
K.\ Kato, C.\ Nakayama and S.\ Usui, 
{\rm Classifying spaces of degenerating mixed Hodge structures},
I., Adv. Stud. Pure Math. {\bf 54} (2009), 187--222, II., Kyoto J. Math. {\bf 51} (2011), 
149--261, III., J. Algebraic Geometry {\bf 22} (2013), 671--772, IV., Kyoto J.\ Math., to appear.

\bibitem{KU}
K.\  Kato and S.\ Usui,
{\rm Classifying spaces of degenerating polarized Hodge structures}, 
Ann.\ Math.\ Studies {\bf 169}, Princeton University Press (2009).

\bibitem{KP}
M.\ Kerr and G.\ Pearlstein, 
{\rm Boundary components of Mumford-Tate domains}, 
Duke Math.\ J.\ {\bf 165-4} (2016), 661--721. 

\bibitem{Ko} T.\ Kohno, 
{\rm Geometry of iterated integrals}, 
Springer contemporary mathematics series {\bf 14} (2012), 295 pages, Maruzen Publ.\ Co., Ltd., (in Japanese). 

\bibitem{Mi}
J.\ S.\ Milne, 
{\rm Canonical models of (mixed) Shimura varieties and automorphic vector bundles}, 
in Automorphic forms, Shimura varieties, and $L$-functions 1 (1990), 283--414.

\bibitem{P}
 D.\ Passman, 
{\rm The algebraic structure of group rings}, 
New York, John Wiley (1977).

\bibitem{U} 
S.\ Usui, 
{\rm Variation of mixed Hodge structure arising from family of logarithmic deformations. II. Classifying space}, 
Duke Math. J. {\bf 51-4} (1984), 851--875.
\end{thebibliography}
\end{document}